\documentclass[11pt]{article}
\usepackage[latin1]{inputenc}
\usepackage{epsfig}
\usepackage{color}
\usepackage[british,english]{babel}
\usepackage{amsthm}
\usepackage{amsmath}
\usepackage{amsfonts}
\usepackage{amssymb}
\usepackage{graphicx}
\def \noame{\noalign{\medskip}}
\newtheorem{corollary}{Corollary}[section]
\newtheorem{definition}[corollary]{Definition}

\newtheorem{lemma}[corollary]{Lemma}
\newtheorem{proposition}[corollary]{Proposition}
\newtheorem{remark}[corollary]{Remark}
\newtheorem{theorem}[corollary]{Theorem}
\newfont{\sBlackboard}{msbm10 scaled 900}

\newcommand{\mylabel}[1]{\label{#1}
            \ifx\undefined\stillediting
            \else \fbox{$#1$}\fi }
\newcommand{\BE}{\begin{equation}}

\newcommand{\EEQ}{\end{equation}}
\newcommand{\rfb}[1]{\mbox{\rm
   (\ref{#1})}\ifx\undefined\stillediting\else:\fbox{$#1$}\fi}

\newfont{\Blackboard}{msbm10 scaled 1200}

\newfont{\roma}{cmr10 scaled 1200}

\def\CC{\rm \hbox{C\kern-.56em\raise.4ex
         \hbox{$\scriptscriptstyle |$}\kern+0.5 em }}

\newcommand{\ep}{\varepsilon}

\def\n{|\kern -.05cm{|}\kern -.05cm{|}}

\newcommand{\mm}    {{\hbox{\hskip 0.5pt}}}

\newcommand{\bluff} {{\hbox{\raise 15pt \hbox{\mm}}}}
\usepackage{fancyhdr}

\makeatletter
\def\section{\@startsection {section}{1}{\z@}{-3.5ex plus -1ex minus
    -.2ex}{2.3ex plus .2ex}{\large\bf}}
\makeatother
\def\be{\begin{equation}}
\def\ee{\end{equation}}

\date{ }
\begin{document}
\thispagestyle{empty}
\title{\bf  Darcy's law for  micropolar fluid flow\\ in a periodic thin porous medium}
\maketitle
\vspace{-50pt}
\author{ \center  Mar\'ia ANGUIANO\\
Departamento de An\'alisis Matem\'atico. Facultad de Matem\'aticas.\\
Universidad de Sevilla, 41012-Sevilla (Spain)\\
anguiano@us.es\\}
\author{ \center  Francisco Javier SU\'AREZ-GRAU\\
Departamento de Ecuaciones Diferenciales y An\'alisis Num\'erico. Facultad de Matem\'aticas.\\
Universidad de Sevilla,  41012-Sevilla (Spain)\\
grau@us.es\\}

\vskip20pt

 \renewcommand{\abstractname} {\bf Abstract}
\begin{abstract} 
In this paper, we extend the Darcy law for micropolar fluid flow in a thin porous medium. This provides a framework for understanding how a fluid's microstructural properties, the geometry of the porous medium and the thickness of the domain (which is significantly smaller than the other dimensions) influence its flow behavior, going beyond the simple pressure-driven flow described by the standard Darcy law.
\end{abstract}
\bigskip\noindent
 {\small \bf AMS classification numbers:}   76A05, 76A20, 76S05, 76M5, 35B27. 

\section {Introduction}\label{S1}
The model of micropolar fluid, proposed by Eringen \cite{Eringen} has been extensively studied both in the
engineering and mathematical literature, due to its practical importance. The
micropolar fluid model describes the motion of numerous real fluids better than the classical
Navier-Stokes equations, because it takes into consideration the microstructure of the fluid particles and capture the effects of its rotation The rotation of the fluid particles is mathematically described by introducing the microrotation field (along with the standard velocity and pressure fields) and, accordingly, a new governing equation coming from the conservation of angular momentum, Lukaszewicz \cite{Luka} (see Section \ref{sec:statement} below for more details). 

The transition between Stokes micropolar equations to Darcy's law assuming Dirichlet conditions on the boundary of the obstacles was rigorously obtained by Aganovi\'c \& Tutek \cite{Aganovic0, Aganovic} and Lukaszewicz \cite[Part III-Chapter 2]{Luka} by means of homogenization (two-scale convergence) techniques, generalizing the classical studies of the derivation of Darcy's laws for Newtonian fluids given in Allaire \cite{Allaire0},  Sanchez-Palencia \cite{Sanchez} and Tartar \cite{Tartar}. We also refer to Bayada {\it et al} \cite{Bayada0Gamouana} where micropolar effects in the coupling of a thin film past a porous medium is studied, and  to Su\'arez-Grau \cite{SG_past} for the study of the derivation of a modified Darcy's law of a micropolar fluid thorugh a porous media with nonzero spin boundary condition on the obstacles.

 On the other hand, the derivation of macroscopic laws for fluids in porous domains with small thickness (the so-called {\it thin porous medium})  is attracting much attention, see for instance Anguiano \& Bunoiu \cite{Ang-Bun2}, Almqvist {\it et al.} \cite{Almqvist}, Anguiano {\it et al.} \cite{Anguiano_Bonnivard_SG, Anguiano_Bonnivard_SG2}, Anguiano \& Su\'arez-Grau \cite{Anguiano_SuarezGrau, Anguiano_SuarezGrau2,  Anguiano_SG_NHM, Anguiano_SG_Lower, Anguiano_SG_sharp}, Fabricius \& Gahn \cite{Fabricius3}, Fabricius {\it et al.} \cite{Fabricius0},   Mei \& Vernescu \cite{Mei}, Prat \&  Aga${\rm \ddot{e}}$sse \cite{Prat},  Su\'arez-Grau \cite{SuarezGrau1}, Yeghiazarian {\it et al.}  \cite{Yeghiazarian} or Zhengan \& Hongxing \cite{Zhengan}.  These studies refer to a different types of Newtonian or non-Newtonian fluid flows through a thin porous medium defined as a bounded perforated 3D domain confined between two parallel plates, where the distance between the plates is very small and the perforation consists of periodically distributed solid cylinders which connect the plates in perpendicular direction.

In this paper, we consider another type of thin porous medium. The domain under consideration is a bounded perforated 3D domain confined between two parallel plates, whose  description 
includes two small parameters:   $\varepsilon$ connected to the microstructure of the domain and $h_\varepsilon$ representing the distance between plates $\varepsilon \ll h_\ep$ (see Figure \ref{fig:domain}). More precisely, we consider
the classical setting of perforated media, i.e. $\varepsilon$-periodically distributed solid (not
connected) obstacles of size $\varepsilon$. This setting has been studied for Newtonian fluids in Bayada {\it et al.} \cite{BayadaThinThin} in a 2D framework and in Su\'arez-Grau \cite{SG_MANA} in a 3D framework. However, as far as the authors know, there is no study in the mathematical literature of micropolar fluids in this type of thin domains.

The goal of this paper is to derive a lower dimensional Darcy's law, adapted for micropolar fluid flow in such a  thin porous medium, which describes the flow of a fluid that not only moves but also rotates, and is influenced by  pressure gradients,  the interaction with the porous medium's structure and the thickness of the domain. This extended Darcy's law incorporates the effects of the fluid's microstructure,  the geometry of the porous medium and the reduction of the dimension, providing a more complete understanding of flow behavior in such complex systems. In the spirit of \cite{BayadaThinThin} and   \cite{SG_MANA}, the key idea is to the derive the sharp a priori estimates using a restriction operator to derive estimates for pressure (Section \ref{sec:estimates}), to use an adaptation of the unfolding method to capture the micro-geometry of the thin porous medium (Section \ref{sec:unfolding}) and to prove the compactness results for the rescaled functions (Section \ref{sec:compactness}).  As a main result formulated in Section \ref{sec:main_thm} (more precisely, see Theorem \ref{mainthm_porous2}), we obtain the homogenized model maintaining at the limit both the effects of velocity and microrotation that represents our main contribution.  
%
%
%
%
%
%
%
%


\section{Geometrical setting and notation} \label{sec:domain} 

\begin{figure}[h!]
\begin{center}
\includegraphics[width=8cm]{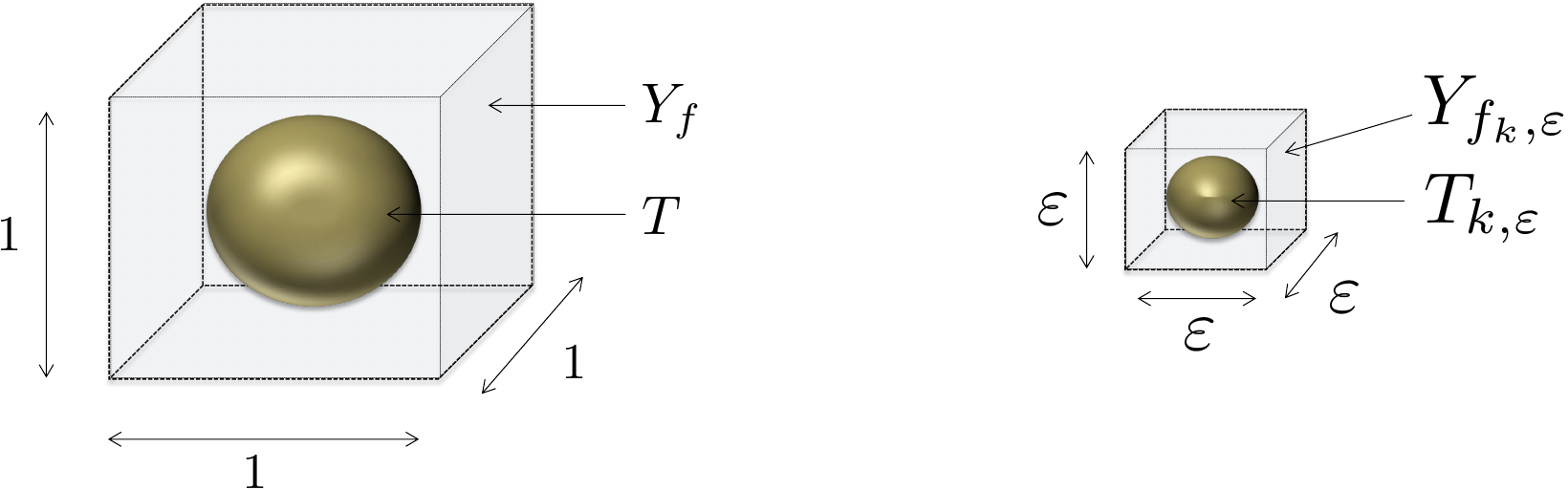}
\end{center}
\vspace{-0.4cm}
\caption{View of the reference cell  $Y$ (left) and the rescaled cell $Y_{k,\ep}$ (right).}
\label{fig:cell}
\end{figure}

\begin{figure}[h!]
\begin{center}
\includegraphics[width=7cm]{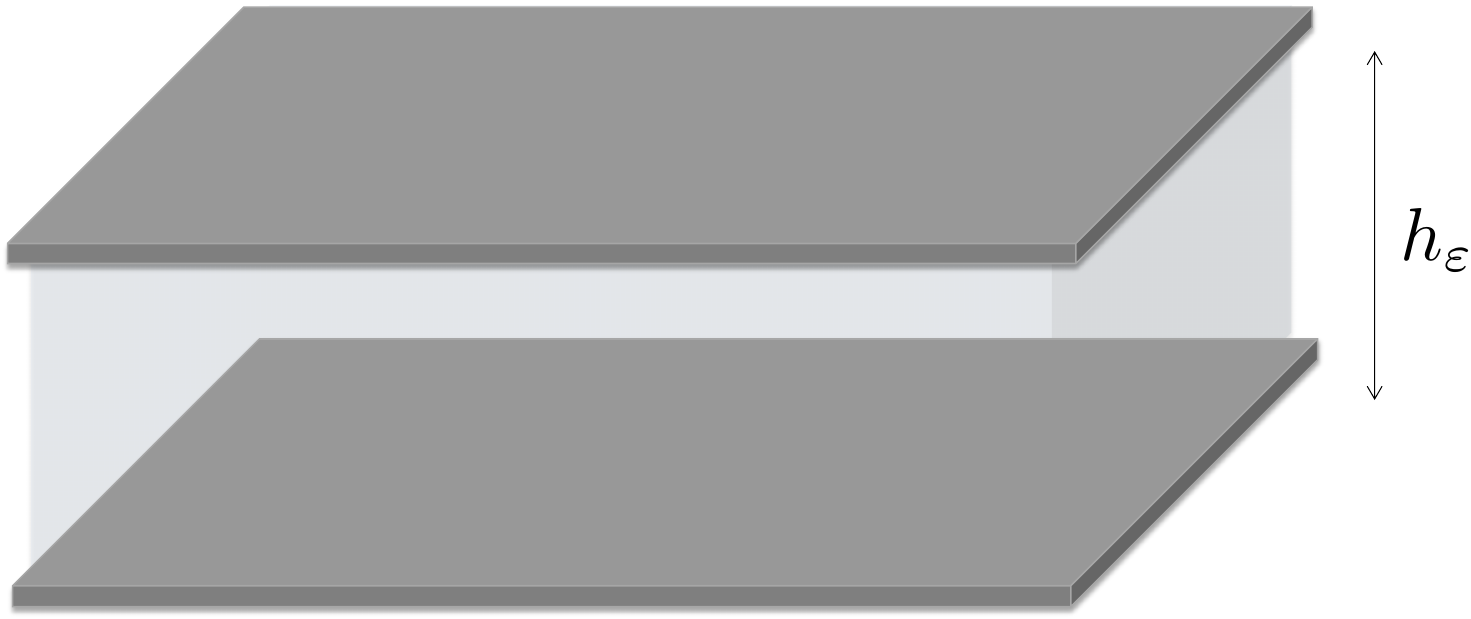}
\hspace{2cm}
\includegraphics[width=7cm]{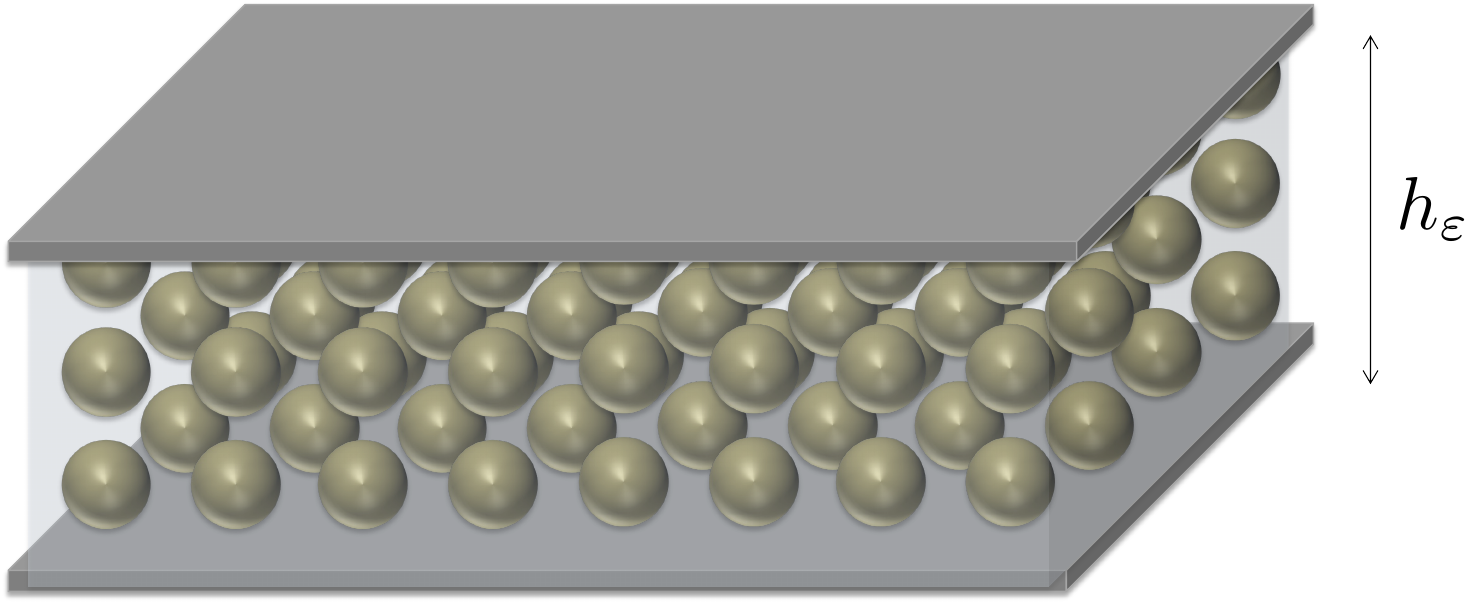}
\end{center}
\caption{View of the domain $Q_\ep$ (left) and $\Omega_{\varepsilon}$ (right).}
\label{fig:domain}
\end{figure}

Let $\omega$ be a smooth, bounded and connected set in $\mathbb{R}^2$. We consider two positive and small parameters $\ep$ and $h_\ep$   (where $h_\ep$ is devoted to tend to zero when $\ep\to0$) satisfying the following relation
\begin{equation}\label{parameters}
\lim_{\ep\to 0} {\ep\over h_\ep}=0,\quad \left(\hbox{i.e.}\quad \ep\ll h_\ep\right).
\end{equation}
We consider   $\Omega_\ep$ is a thin porous medium, where to be described,  we consider the parameters $\ep$ and $h_\ep$ satisfying (\ref{parameters}). Thus, we consider a thin layer of height $h_\ep$ which is perforated by $\ep$-periodic distributed obstacles of size $\ep$. The thin layer without microstructure is denoted by $Q_\ep$ (see Figure \ref{fig:domain}), i.e.
\begin{equation}\label{Qep}
Q_\ep=\omega\times (0,h_\ep).
\end{equation} 
Let us now give a better description of the microstructure of the thin layer. We denote $Y=(-1/2,1/2)^3$ the unitary cube in $\mathbb{R}^3$ as the reference cell and  $T$ an open connected subset of $Y$ with a smooth boundary  $\partial T$ such that $\overline T\subset Y$. We denote $Y_f=Y\setminus \overline T$. Thus, for $k\in\mathbb{Z}^3$, each cell $Y_{k,\varepsilon}=\varepsilon k+\varepsilon Y$ is similar to the unit cell $Y$ rescaled to size $\varepsilon$ and $T_{k,\varepsilon}=\varepsilon k+\varepsilon T$ is similar to $T$ rescaled to size $\varepsilon$. We denote $Y_{f_k,\varepsilon}=Y_{k,\varepsilon}\setminus \overline T_{k,\varepsilon}$ (see Figure \ref{fig:cell}).

We denote by $\tau(\overline T_{k,\ep})$ the set of all translated images of $\overline T_{k,\ep}$, i.e. the set $\tau(\overline T_{k,\ep})$ represents the obstacles in $\mathbb{R}^3$.   The thin porous media $\Omega_\ep$ is defined by  (see Figure \ref{fig:domain})
\begin{equation}\label{Omegaep}
\Omega_\ep=Q_\ep\setminus \bigcup_{k\in \mathcal{K}_\ep}\overline T_{k,\varepsilon},
\end{equation} where   $\mathcal{K}_\ep:=\left\{k\in\mathbb{Z}^3\,:\, Y_{k,\varepsilon}\cap Q_\ep\neq \emptyset\right\}$. By construction,  $\Omega_\varepsilon$ is a periodically perforated channel with obstacles of the same size as the period.
  We make the assumption that the obstacles $\tau(\overline T_{k,\ep})$ do no intersect the boundary $\partial Q_\ep$. We denote by $T_{\varepsilon}$ the set of all the obstacles contained in $\Omega_\ep$. Then,   $T_\ep$ is a finite union of obstacles, i.e. 
 $$T_{\varepsilon}=\bigcup_{k\in \mathcal{K}_\ep}\overline T_{k,\ep}.$$

\begin{figure}[h!]
\begin{center}
\includegraphics[width=4cm]{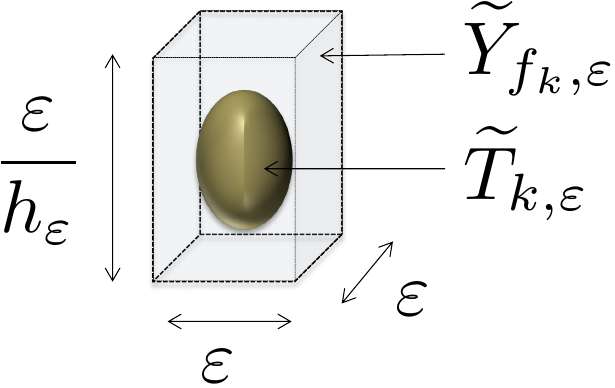}
\hspace{2.5cm}
{\includegraphics[width=7cm]{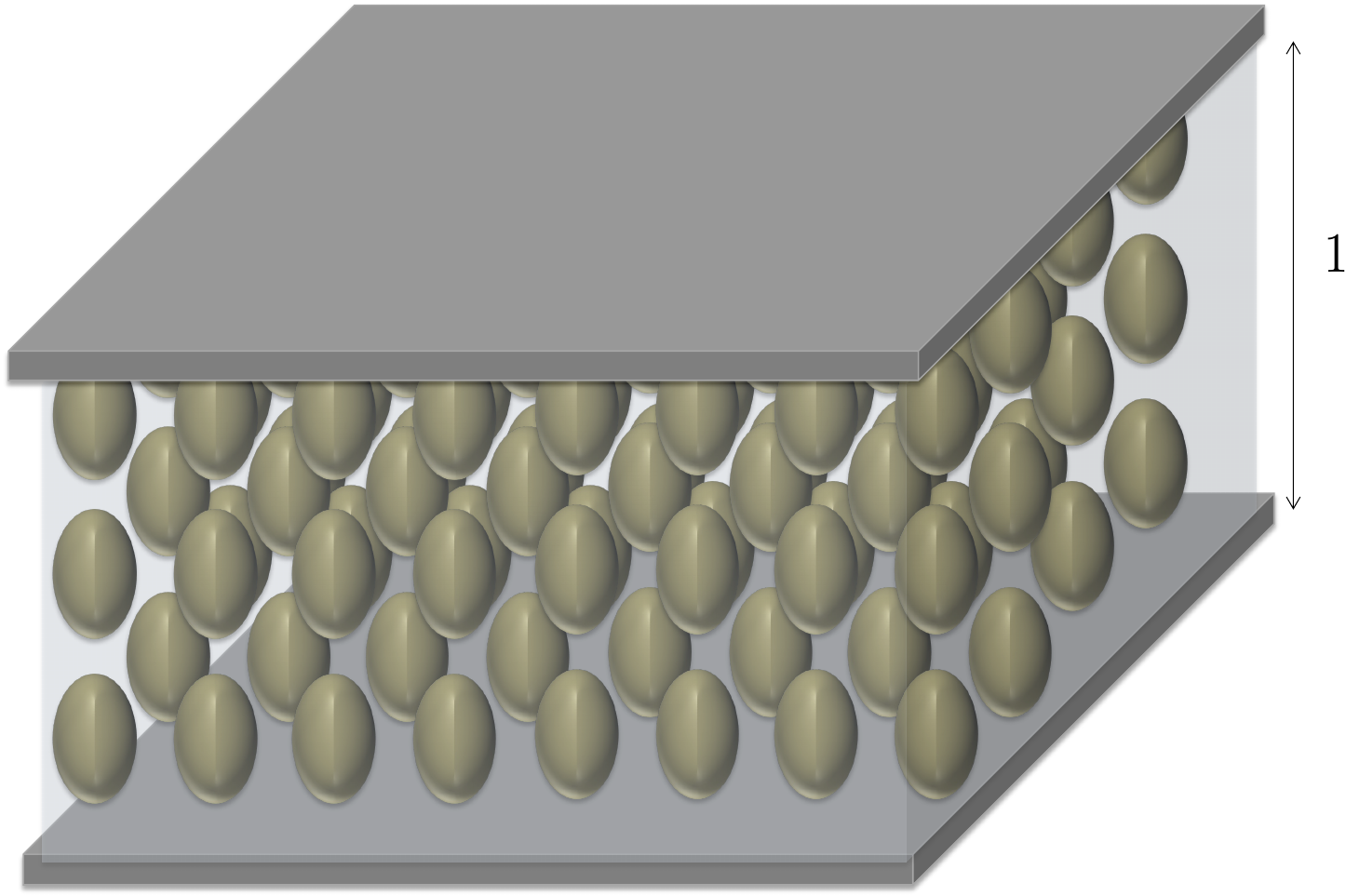}}
\end{center}
\vspace{-0.4cm}
\caption{View of the rescaled reference cell  $\widetilde Y_{k,\ep}$ (left) and  the rescaled domain $\widetilde\Omega_\ep$ (right).}
\label{fig:dilated}
\end{figure}

\noindent As usual when we deal with thin domains, we will use the dilatation in the variable $x_3$ given by
\begin{equation}\label{dilatacion_p}
z'= x',\quad z_3={x_3\over h_\ep},\quad \forall\,x\in \Omega_\ep.
\end{equation}
Then, we define the rescaled porous media  $\widetilde \Omega_\ep$ by   (see Figure \ref{fig:dilated})
\begin{equation}\label{Omega_tilde}
\widetilde\Omega_\ep=\left\{z=(z',z_3)\in \mathbb{R}^2\times \mathbb{R}\,:\, (z',h_\ep z_3)\in \Omega_\ep\right\}.
\end{equation}
We also introduce the rescaled sets $\widetilde Y_{k,\ep}$ by  (see Figure \ref{fig:dilated})
$$\widetilde Y_{{k},\ep}=\left\{y\in\mathbb{R}^2\times \mathbb{R}\,:\, (y',h_\ep y_3)\in Y_{k,\varepsilon}\right\},$$
and, in the same way, we define the rescaled fluid part $\widetilde Y_{f_{k},\ep}$, the rescaled solid  part $\widetilde T_{k,\varepsilon}$ of $\widetilde Y_{{k},\ep}$ and the union of rescaled obstacles $\widetilde T_\ep$.

Finally,  we denote by $\Omega$ the domain with fixed height without microstructure, i.e.
$$
\Omega=\omega\times (0,1),
$$
and by $\Gamma_0$ and $\Gamma_1$ the bottom and top boundaries of $\Omega$, i.e. 
$$\Gamma_0=\omega\times \{0\},\quad \Gamma_1=\omega\times \{1\}.$$
 
 Along this paper, the points $x\in \mathbb{R}^3$ will be decomposed as $x=(x',x_3)$ with $x'\in \mathbb{R}^2$, $x_3\in \mathbb{R}$. We also use the notation $x'$ to denote a generic vector of $\mathbb{R}^2 $. Let us consider a vectorial function $\varphi=(\varphi', \varphi_3)$, $\varphi'=(\varphi_1, \varphi_2)$ and a scalar function $\psi$.    We use  the following operators 
$$\Delta\varphi=\Delta_{x'}\varphi+\partial_{x_3}^2\varphi,\quad {\rm div}(\varphi)={\rm div}_{x'}(\varphi')+\partial_{x_3}\varphi_3,\quad \nabla\psi=(\nabla_{x'}\psi, \partial_{x_3}\psi)^t. $$
$$
{\rm rot}(\varphi)=(\partial_{x_2}\varphi_3-\partial_{x_3}\varphi_2, \partial_{x_3}\varphi_1-\partial_{x_1}\varphi_3, \partial_{x_1}\varphi_2-\partial_{x_3}\varphi_1)^t. $$
 For a vectorial function $\widetilde \varphi=(\widetilde \varphi',\widetilde \varphi_3)$ and a scalar function $\widetilde \psi$ obtained respectively from $\varphi$ and $\psi$ by using the rescaling (\ref{dilatacion_p}) in the set $\Omega_\ep$, we will denote   $\Delta_{h_\ep}$, $D_{h_\ep}$, $\nabla_{h_\ep}$,  ${\rm div}_{h_\ep}$ and ${\rm rot}_{h_\ep}$ as follows
$$\begin{array}{c}
\displaystyle (D_{ h_\ep}\widetilde \varphi)_{i,j}=\partial_{z_j}\widetilde \varphi_i,\quad i=1, 2, 3,\ j=1, 2, \displaystyle (D_{ h_\ep}\widetilde \varphi)_{i,3}=h_\ep^{-1}\partial_{z_3}\widetilde \varphi_i\hbox{ for }i=1,2, 3,\\
\noame
\displaystyle \Delta_{h_\ep}\widetilde\varphi=\Delta_{z'}\widetilde\varphi+h_\ep^{-2}\partial_{z_3}^2\widetilde\varphi,\quad \displaystyle  {\rm div}_{h_\ep} (\widetilde\varphi)={\rm div}_{z'}(\widetilde\varphi')+h_\ep^{-1}\partial_{z_3}\widetilde\varphi_3,\quad
\nabla_{ h_\ep}\widetilde\psi=(\nabla_{z'}\widetilde\psi,h_\ep^{-1}\partial_{z_3}\widetilde\psi)^t,\\
\noame
\displaystyle
{\rm rot}_{h_\ep}(\widetilde\varphi)=({\rm rot}_{z'}(\widetilde \varphi_3) +h_\ep^{-1}{\rm rot}_{z_3}(\widetilde \varphi'),{\rm Rot}_{z'}(\widetilde \varphi'))^t ,
\end{array}$$
where, denoting $(\widetilde \varphi')^\perp=(-\widetilde \varphi_2,\widetilde \varphi_1)^t$, we define
$${\rm rot}_{x'}(\widetilde \varphi_3)=(\partial_{z_2}\widetilde \varphi_3, -\partial_{z_1}\widetilde \varphi_3)^t,\quad {\rm rot}_{z_3}(\widetilde \varphi')=(\partial_{z_3}\widetilde \varphi')^\perp,\quad {\rm Rot}_{z'}(\widetilde \varphi')=\partial_{z_1}\widetilde \varphi_2-\partial_{z_2}\widetilde \varphi_1.$$

 Let $C^\infty_{\rm per}(Y)$ be the space of infinitely differentiable functions in $\mathbb{R}^3$ that are $Y$-periodic. By $L^2_{\rm per}(Y)$ (resp. $H^1_{\rm per}(Y)$) we denote its completion in the norm $L^2(Y)$ (resp. $H^1(Y)$). We denote by $L^{2}_0$ the space of functions of $L^{2}$ with null integral and by $L^{2}_{0,{\rm per}}(Y)$  the space of functions in $L^{2}_{\rm per}(Y)$ with zero mean value.

We denote by $:$ the full contraction of two matrices, i.e. for $A=(a_{ij})_{1\leq i,j\leq 3}$ and $B=(a_{ij})_{1\leq i,j\leq 3}$, we have $A:B=\sum_{i,j=1}^3a_{ij}b_{ij}$.  The canonical basis in $\mathbb{R}^3$ is denoted by $\{e_1,e_2, e_3\}$.

Finally, we denote by $O_\ep$ a generic real sequence, which tends to zero with $\ep$ and can change from line to line, and by $C$ a generic positive constant which also can change from line to line.

\section{Statement of the problem}\label{sec:statement}  In view of the application we want to model, we can assume a small Reynolds number and neglect the inertial terms in the governing equations. We study the motion of an incompressible micropolar fluid in the whole fluid domain $\Omega_\ep$ governed by the dimensionless  linearized micropolar Stokes law  
\begin{equation}\label{system_1introad}
\left\{\begin{array}{rl}
- \Delta u_\ep +\nabla p_\ep =2N^2  {\rm rot}( w_\ep )+  f_\ep,  &\\
\\
{\rm div}(   u_\ep)=0,&\\
\\
\displaystyle -R_M\Delta   w_\ep + 4N^2 w_\ep =2N^2 {\rm rot}( u_\ep)+  g_\ep, & 
\end{array}\right.
\end{equation}
where the subscript $\ep$ is added to the unknowns to stress the dependence of the solution on the small parameter, with no-slip and no-spin boundary boundary conditions
\begin{equation}\label{system_1intro2ad}
  u_\ep=0,\quad   w_\ep=0,\quad \hbox{on }\partial Q_\ep\cup \partial T_\ep,
\end{equation}
where $u_\ep$, $p_\ep$ and $w_\ep$ are respectively the velocity, the pressure and the angular rotation velocity of the particles in the fluid (called the microrotation).     $N^2$ is called a coupling number and characterizes the coupling of equations (\ref{system_1introad})$_1$ and (\ref{system_1introad})$_3$, and it is of order $\mathcal{O}(1)$ such that  $0<N^2<1$. The second dimensionless parameter, denoted by $R_M$ is, in fact, related to the characteristic length of the microrotation effects and will be compared with small parameter $\ep$ (see condition (\ref{RM}) below).\\


Our aim is to describe the asymptotic behavior of the velocity $u_\ep$, the pressure
$p_\ep$ and microrotation $w_\ep$  as $\ep$ tends to zero and identify an homogenized model coupling the effects of the thickness and the microgeometry of the domain.  To do this,  we will use the equivalent weak variational formulation of (\ref{system_1introad})--(\ref{system_1intro2ad}), which is the following one:

Find $u_\ep\in H^1_0(\Omega_\ep)^3$, $p_\ep\in L^2_0(\Omega_\ep)$ and $w_\ep\in H^1_0(\Omega_\ep)^3$  such that
\begin{equation}\label{form_var}
\left\{\begin{array}{l}
\displaystyle\int_{\Omega_\ep}D u_\ep: D\varphi\,dx-\int_{\Omega_\ep}p_\ep\,{\rm div}(\varphi)\,dx=2N^2\int_{\Omega_\ep}{\rm rot}(w_\ep)\cdot \varphi\,dx+\int_{\Omega_\ep}f_\ep\cdot \varphi\,dx,\\
\noame
\displaystyle\int_{\Omega_\ep}{\rm div}(u_\ep)\zeta\,dx=0,\\
\noame
\displaystyle R_M\int_{\Omega_\ep}D w_\ep: D\psi\,dx+4N^2\int_{\Omega_\ep}w_\ep\cdot\psi\,dx=2N^2\int_{\Omega_\ep}{\rm rot}(u_\ep)\cdot \psi\,dx+\int_{\Omega_\ep}g_\ep\cdot \psi\,dx,\\
\end{array}
\right.
\end{equation}
\indent for every $\varphi\in H^1_0(\Omega_\ep)^3$, $\zeta\in L^2(\Omega_\ep)$ and $\psi\in H^1_0(\Omega_\ep)^3$.\\

It is well known (see for instance  \cite{Luka}) that, for every $\ep > 0$, problem (\ref{system_1introad})-(\ref{system_1intro2ad}) has a unique weak solution $(u_\ep,p_\ep,w_\ep)\in H^1_0(\Omega_\ep)^3\times L^2_0(\Omega_\ep)\times H^1_0(\Omega_\ep)^3$. 
\\

We remark that different asymptotic behaviors of the flow can be deduced depending on the order of magnitude of the dimensionless parameter $\ep$. Indeed, if we compare the characteristic number $R_M$ with small parameter $\ep$, three different asymptotic situations can be formally identified (see e.g. \cite{Bayada96, Bayada01}). We consider the most interesting one, which leads to a strong coupling at main order, namely the regime
\begin{equation}\label{RM}
R_M=\ep^2 R_c,\quad \hbox{with}\quad R_c=\mathcal{O}(1).
\end{equation}
 Moreover, we assume  
\begin{equation}\label{fep}
 f_\ep=(f'(x'),0)\quad\hbox{with}\quad f'\in L^{2}(\omega)^2, \quad   g_\ep= \ep (g'(x'),0),\quad \hbox{with}\quad g'\in L^2(\omega)^2,
\end{equation}
which is usual when we deal with thin domains. Since the thickness of the domain is small, then the vertical component of the force can be neglected, and moreover, the force can be considered independent of the vertical variable. 


In order to find the limit problem when $\ep$ tends to zero, we introduce the rescaling given by (\ref{dilatacion_p}). Using this rescaling, we can define  $\widetilde u_\ep \in H^1_0(\widetilde \Omega_\ep)^3$, $\widetilde p_\ep\in L^{2}_0(\widetilde \Omega_\ep)$ and $\widetilde w_\ep \in H^1_0(\widetilde \Omega_\ep)^3$ by 
\begin{equation}\label{unk_dilat}
\displaystyle \widetilde u_\ep(z)=u_\ep(z', h_\ep z_3),\quad  \widetilde p_\ep(z)=p_\ep(z', h_\ep z_3),\quad \widetilde w_\ep(z)=w_\ep(z', h_\ep z_3)\quad\hbox{if }z\in \widetilde\Omega_\ep\,,
\end{equation}
so the rescaled weak variational formulation is the following: 

Find $\widetilde u_\ep\in H^1_0(\widetilde\Omega_\ep)^3$, $\widetilde p_\ep\in L^2_0(\widetilde\Omega_\ep)$ and $\widetilde w_\ep\in H^1_0(\widetilde\Omega_\ep)^3$  such that
\begin{equation}\label{form_var_tilde}
\left\{\begin{array}{l}
\displaystyle\int_{\widetilde\Omega_\ep} D_{h_\ep} \widetilde u_\ep: D_{h_\ep}\widetilde \varphi\,dz-\int_{\widetilde\Omega_\ep} \widetilde p_\ep\,{\rm div}_{h_\ep}(\widetilde \varphi)\,dz
=2N^2\int_{\widetilde\Omega_\ep} {\rm rot}_{h_\ep}(\widetilde w_\ep)\cdot \widetilde \varphi\,dz +\int_{\widetilde\Omega_\ep}  f'\cdot \widetilde \varphi'\,dz,\\
\\
\displaystyle\int_{\widetilde\Omega_\ep} {\rm div}_{h_\ep}(\widetilde u_\ep)\widetilde \zeta\,dz =0,\\
\\
\displaystyle \ep^2 R_c\int_{\widetilde\Omega_\ep}  D_{h_\ep}\widetilde w_\ep: D_{h_\ep}\widetilde\psi\,dz 
+4N^2\int_{\widetilde\Omega_\ep}  \widetilde w_\ep\cdot\widetilde\psi\,dz=2N^2\int_{\widetilde\Omega_\ep} {\rm rot}_{h_\ep}(\widetilde u_\ep)\cdot\widetilde \psi\,dz +\ep \int_{\widetilde\Omega_\ep}g'\cdot\widetilde \psi'\,dz,
\end{array}\right.
\end{equation}
\indent for every $\widetilde \varphi\in H^1_0(\widetilde \Omega_\ep)^2$, $\widetilde \zeta\in L^2(\widetilde \Omega_\ep)$ and $\widetilde \psi\in H^1_0(\widetilde \Omega_\ep)$.\\

\section{A priori estimates}\label{sec:estimates}
Let us begin with the  Poincar\'e   inequality in the thin porous medium $\Omega_\ep$ (see \cite{Tartar} for classical version) and the relation between estimates of rotational and derivative.
\begin{lemma}[Lemma 5.1 in \cite{SG_MANA}] \label{Lemma_Poincare} For every $\varphi\in H^1_0(\Omega_{\varepsilon})^3$ there exists a positive constant $C_P$, independent of $\ep$, such that, 
\begin{equation}\label{Poincare1}
\|\varphi\|_{L^2(\Omega_\ep)^3}\leq C_P\ep\|D \varphi\|_{L^2(\Omega_\ep)^{3\times 3}}.\end{equation}
\end{lemma}
 
\begin{lemma}[\cite{Luka}] \label{Lemma_rot_D} For every $\varphi\in H^1_0(\Omega_{\varepsilon})^3$ then it holds
\begin{equation}\label{rot1}
\|{\rm rot}(\varphi)\|_{L^2(\Omega_\ep)^3}\leq  \|D \varphi\|_{L^2(\Omega_\ep)^{3\times 3}}.\end{equation}
Moreover, if ${\rm div}(\varphi)=0$ in $\Omega_\ep$, then
\begin{equation}\label{rot2}
\|{\rm rot}(\varphi)\|_{L^2(\Omega_\ep)^3}= \|D \varphi\|_{L^2(\Omega_\ep)^{3\times 3}}.\end{equation}
\end{lemma}

\begin{lemma}[Estimates for velocity and microrotation]\label{Lem:estimates_velocity}There exists a constant $C>0$ independent of $\ep$, such that if $(u_\ep,w_\ep)\in H^1_0(\Omega_\ep)^2\times H^1_0(\Omega_\ep)$ is the solution of problem (\ref{system_1introad})-(\ref{system_1intro2ad}), it holds
\begin{equation}\label{estim_Du}
\|u_\ep\|_{L^2(\Omega_\ep)^{3}}\leq C\ep^2 h_\ep^{1\over 2},\quad 
\|Du_\ep\|_{L^2(\Omega_\ep)^{3\times 3}}\leq C\ep h_\ep^{1\over 2},
\end{equation}
\begin{equation}\label{estim_Dw}
\|w_\ep\|_{L^2(\Omega_\ep)^{3}}\leq C \ep h_\ep^{1\over 2},\quad 
\|Dw_\ep\|_{L^2(\Omega_\ep)^{3\times 3}}\leq C h_\ep^{1\over 2}.
\end{equation}
Moreover, applying the dilatation (\ref{dilatacion_p}), we get
\begin{equation}\label{estim_Du_tilde}
\|\widetilde u_\ep\|_{L^2(\widetilde\Omega_\ep)^{3}}\leq C\ep^2,\quad 
\|D_{h_\ep} \widetilde u_\ep\|_{L^2(\widetilde\Omega_\ep)^{3\times 3}}\leq C\ep,
\end{equation}
\begin{equation}\label{estim_Dw_tilde}
\|\widetilde w_\ep\|_{L^2(\widetilde \Omega_\ep)^{3}}\leq C\ep,\quad 
\|D_{h_\ep} \widetilde w_\ep\|_{L^2(\widetilde\Omega_\ep)^{3\times 3}}\leq C.
\end{equation}
\end{lemma}

\begin{proof} We divide the proof in two steps. In the first step, we derive estimates for velocity and then, for microrotation.\\

{\it Step 1. Velocity estimates}. Using $u_\ep$ as test function in (\ref{form_var})$_1$, we have
\begin{equation}\label{estim_proof_1}
 \|D u_\ep\|_{L^2(\Omega_\ep)^{3\times 3}}^2=2N^2\int_{\Omega_\ep}{\rm rot}(w_\ep)\cdot u_\ep\,dx+\int_{\Omega_\ep}f'\cdot u_\ep'\,dx.
\end{equation}
Using $\int_{\Omega_\ep}{\rm rot}(w_\ep)\cdot u_\ep\,dx=\int_{\Omega_\ep}{\rm rot}(u_\ep)\cdot w_\ep\,dx$  (see \cite{Luka}),    the Cauchy-Schwarz inequality, equality (\ref{rot2}) and the assumption of $f'$ given in (\ref{fep}), we obtain 
$$\begin{array}{rl}
\displaystyle\|D u_\ep\|_{L^2(\Omega_\ep)^{3\times 3}}^2\leq &\displaystyle 2N^2\int_{\Omega_\ep}w_\ep\cdot {\rm rot}(u_\ep)\,dx+\int_{\Omega_\ep}f'\cdot u_{\ep}'\,dx\\
\noame
=&\displaystyle
2N^2\|w_\ep\|_{L^2(\Omega_\ep)^3}\|D u_\ep\|_{L^2(\Omega_\ep)^{3\times 3}}+ C_P h_\ep^{1\over 2}\ep \|f'\|_{L^2(\omega)^2}\|D u_\ep\|_{L^2(\Omega_\ep)^{3\times 3}}.
\end{array}$$
Therefore,  we get
\begin{equation}\label{estim_proof_3_vel}
\|D u_\ep\|_{L^2(\Omega_\ep)^{3\times 3}}\leq 2N^2\|w_\ep\|_{L^2(\Omega_\ep)^3} + C_P h_\ep^{1\over 2}\ep \|f'\|_{L^2(\omega)^2}.
\end{equation}
On the other hand, taking $w_\ep$ as test function in (\ref{form_var})$_3$, we have
\begin{equation}\label{estim_proof_2}
 \ep^2 R_c \|D w_\ep\|_{L^2(\Omega_\ep)^{3\times 3}}^2+4N^2\|w_\ep\|_{L^2(\Omega_\ep)^3}=2N^2\int_{\Omega_\ep}{\rm rot}(u_\ep)\cdot w_\ep\,dx+\ep\int_{\Omega_\ep}g'\cdot w_\ep'\,dx.
\end{equation}
Applying the Cauchy-Schwarz inequality  and the assumption of $g'$ given in (\ref{fep}), we get
$$
\begin{array}{rl}
2N^2\|w_\ep\|_{L^2(\Omega_\ep)^3}^2\leq &\displaystyle N^2\int_{\Omega_\ep}{\rm rot}(u_\ep)\cdot  w_\ep\,dx+{\ep\over 2}\int_{\Omega_\ep}g'\cdot w_\ep'\,dx\\
\noame
\leq &\displaystyle N^2\|w_\ep\|_{L^2(\Omega_\ep)^3}\|D u_\ep\|_{L^2(\Omega_\ep)^{3\times 3}}+{1\over 2}\ep h_\ep^{1\over 2}\|g'\|_{L^2(\omega)^2}\|w_\ep\|_{L^2(\Omega_\ep)^3},
\end{array}
$$
which implies
\begin{equation}\label{estim_proof_3}
2N^2\|w_\ep\|_{L^2(\Omega_\ep)^3}^2\leq  N^2 \|D u_\ep\|_{L^2(\Omega_\ep)^{3\times 3}}+{1\over 2}\ep h_\ep^{1\over 2}\|g'\|_{L^2(\omega)^2}.
\end{equation}
Then, putting (\ref{estim_proof_3}) into (\ref{estim_proof_3_vel}), we deduce
$$(1-N^2)\|D u_\ep\|_{L^2(\Omega_\ep)^{2\times 2}}\leq \ep h_\ep^{1\over 2}\left(C_P \|f'\|_{L^2(\omega)^2}+{1\over 2}\|g'\|_{L^2(\omega)^2}  \right),$$
which, since $1-N^2>0$,   gives (\ref{estim_Du})$_2$ and then, by using (\ref{Poincare1}), we deduce  (\ref{estim_Du})$_1$.
\\

{\it Step 2. Microrotation estimates.} From (\ref{estim_proof_2}), using   $\int_{\Omega_\ep}{\rm rot}(u_\ep)\cdot w_\ep\,dx=\int_{\Omega_\ep}{\rm rot}(w_\ep)\cdot u_\ep\,dx$,    the Cauchy-Schwarz inequality, estimate (\ref{rot2}) and the assumption of $g'$ given in (\ref{fep}), we obtain 
$$
\begin{array}{rl}
\ep^2 R_c \|Dw_\ep\|_{L^2(\Omega_\ep)^{3\times 3}}^2\leq &\displaystyle 2 N^2\|u_\ep\|_{L^2(\Omega_\ep)^3}\|D w_\ep\|_{L^2(\Omega_\ep)^{3\times 3}}+ C_P\ep^2 h_\ep^{1\over 2}\|g'\|_{L^2(\omega)^2}\|Dw_\ep\|_{L^2(\Omega_\ep)^{3\times 3}},
\end{array}
$$
which, by using estimates of $u_\ep$ given in (\ref{estim_Du})$_1$, gives
$$
\begin{array}{rl}
\ep^2 R_c \|Dw_\ep\|_{L^2(\Omega_\ep)^{3\times 3}}\leq &\displaystyle 2 N^2\ep^2 h_\ep^{1\over 2} + C_P\ep^2 h_\ep^{1\over 2}\|g'\|_{L^2(\omega)^2}.
\end{array}
$$
This implies (\ref{estim_Dw})$_2$ and by using (\ref{Poincare1}), we get (\ref{estim_Dw})$_1$.

{\it Step 3. Estimates for dilated velocity and microrotation.}  Estimates for rescaled velocity  and microrotation  (\ref{estim_Du_tilde})-(\ref{estim_Dw_tilde})  in $\widetilde \Omega_\ep$ are obtained directly from  the estimates (\ref{estim_Du})-(\ref{estim_Dw}) by applying the change of variables (\ref{dilatacion_p}), just taking into account that
$$\begin{array}{c}
\displaystyle \|\varphi_\ep\|_{L^2(\Omega_\ep)^3}=h_\ep^{{1\over 2}}\|\widetilde \varphi_\ep\|_{L^2(\widetilde \Omega_\ep)^3},\quad   \|D  \varphi_\ep\|_{L^2(  \Omega_\ep)^{3\times 3}}=h_\ep^{{1\over 2}}\|D_{h_\ep} \widetilde \varphi_\ep\|_{L^2(\widetilde \Omega_\ep)^{3\times 3}},
\end{array}$$
for $\varphi=\{u, w\}.$
\end{proof}

  Next, we derive a priori estimates for the pressure in the porous part. To do this, we need to extend the pressure to the whole thin film $Q_\ep$ (which also depends on $\ep$).  To do this, we  recall a result from   \cite[Lemma 3.3]{BayadaThinThin} and  \cite[Lemma 5.3]{SG_MANA} (see also \cite[Lemma 3.5.]{Anguiano_SG_fissure2025} for a generalization for $W^{1,q}$, $1<q<+\infty$), which introduces a restriction operator $\mathcal{R}^\ep$ from $H^1_0(Q_\ep)^3$ into $H^1_0(\Omega_\ep)^3$.
 \begin{lemma}[\cite{BayadaThinThin, SG_MANA}] \label{restriction_operator}
There exists  a (restriction) operator $\mathcal{R}^\ep$ acting from $H^1_0(Q_\ep)^3$ into $H^1_0(\Omega_\ep)^3$ such that
\begin{enumerate}
\item $\mathcal{R}^\ep \varphi=\varphi$, if $\varphi \in H^1_0(\Omega_\ep)^3$.
\item ${\rm div}(\mathcal{R}^\ep \varphi)=0\hbox{  in }\Omega_\ep$, if ${\rm div}(\varphi)=0\hbox{  on }Q_\ep$.
\item For every $\varphi\in H^1_0(Q_\ep)^3$, there exists a positive constant $C$, independent of $\varphi$ and $\ep$, such that
\begin{equation}\label{estim_restricted}
\begin{array}{l}
\|\mathcal{R}^\ep \varphi\|_{L^2(\Omega_\ep)^{3}}+ \ep\|D \mathcal{R}^\ep \varphi\|_{L^2(\Omega_\ep)^{3\times 3}} \leq C\left(\|\varphi\|_{L^2(Q_\ep)^3}+\ep \|D\varphi\|_{L^2(Q_\ep)^{3\times 3}}\right)\,.
\end{array}
\end{equation}
\end{enumerate}
\end{lemma}
 
By applying the change of variables (\ref{dilatacion_p}) to estimates (\ref{estim_restricted}), we directly obtain the following result.

\begin{lemma}[\cite{BayadaThinThin, SG_MANA}] \label{restriction_operator2} Setting $\mathcal{\widetilde R}^\ep\widetilde\varphi=\mathcal{R}^\ep\varphi$ for any $\widetilde \varphi\in H^1_0(\Omega)^2$, where $\widetilde \varphi$ is obtained from $ \varphi$ by using the change of variables (\ref{dilatacion_p}),  and $\mathcal{R}^\ep$ is defined in Lemma \ref{restriction_operator}, we have the following estimates:
\begin{equation}\label{estim_restricted_tilde}
\begin{array}{l}
\|\mathcal{\widetilde R}^\ep\widetilde  \varphi\|_{L^2(\widetilde \Omega_\ep)^{3}} \leq \displaystyle C\left(\|\widetilde \varphi\|_{L^2(\Omega)^3}+\ep \|D_{z'}\widetilde \varphi\|_{L^2(\Omega)^{3\times 2}}+{\ep\over h_\ep} \|\partial_{z_3}\widetilde \varphi\|_{L^2(\Omega)^{3}}\right)\,,\\
\noame
 \|D_{z'} \mathcal{\widetilde R}^\ep\widetilde  \varphi\|_{L^2(\widetilde \Omega_\ep)^{3\times 2}} \leq C\left(\ep^{-1}\|\widetilde \varphi\|_{L^2(\Omega)^3}+  \|D_{z'}\widetilde \varphi\|_{L^2(\Omega)^{3\times 2}}+ h_\ep^{-1} \|\partial_{z_3}\widetilde \varphi\|_{L^2(\Omega)^{3}}\right)\,,\\
\noame
 \|\partial_{z_3} \mathcal{R}^\ep \varphi\|_{L^2(\widetilde\Omega_\ep)^{3\times 3}} \leq \displaystyle C\left({h_\ep\over \ep}\|\widetilde \varphi\|_{L^2(\Omega)^3}+  h_\ep\|D_{z'}\widetilde \varphi\|_{L^2(\Omega)^{3\times 2}}+  \|\partial_{z_3}\widetilde \varphi\|_{L^2(\Omega)^{3}}\right)\,,
\end{array}
\end{equation}
which, from relation (\ref{parameters}), imply
\begin{equation}\label{estim_restricted2}
\begin{array}{l}
\displaystyle 
\|\mathcal{\widetilde R}^\ep \widetilde \varphi\|_{L^2(\widetilde \Omega_\ep)^{3}}\leq C\|\widetilde \varphi\|_{H^1_0(\Omega)^3},\quad \|D_{h_\ep} \mathcal{\widetilde R}^\ep \widetilde \varphi\|_{L^2(\widetilde \Omega_\ep)^{3\times 3}} \leq C\ep^{-1}\|\widetilde \varphi\|_{H^1_0(\Omega)^3}.
\end{array}
\end{equation}

\end{lemma}
 
By using the restriction operator,  we give the existence of an extended pressure to $Q_\ep$ by duality arguments, and also derive the corresponding estimates.

\begin{lemma} \label{Estimates_extended_lemma}  There exists an extension $P_\ep\in L^{2}_0(Q_\ep)$ of the pressure $p_\ep$. Moreover, defining the rescaled and extended pressure $\widetilde P_\ep\in L^{2}_0(\Omega)$ obtained from $P_\ep$ by using the change of variables (\ref{dilatacion_p}), then 
there exists a positive constant $C$ independent of $\ep$, such that 
\begin{equation}\label{esti_P}
\|\widetilde P_\ep\|_{L^{2}(\Omega)}\leq C ,\quad 
\|\nabla_{h_\ep}  \widetilde P_\ep\|_{H^{-1}(\Omega)^3}\leq C .
\end{equation}
\end{lemma}
\begin{proof}  We divide the proof in two steps. In the first step, we extend the pressure $p_\ep$ and next, we derive estimates for the rescaled and extended pressure.

{\it Step 1. Extension of $p_\ep$ to $Q_\ep$}. Using the restriction operator $\mathcal{R}^\ep$ given in  Lemma \ref{restriction_operator}, we   introduce $F_\ep$ in $H^{-1}(Q_\ep)^3$ in the following way
\begin{equation}\label{F}\langle F_\varepsilon, \varphi\rangle_{H^{-1}(Q_\varepsilon)^3, H^1_0(Q_\ep)^3}=\langle \nabla p_\varepsilon, \mathcal{R}^\varepsilon \varphi\rangle_{{H^{-1}(\Omega_\varepsilon)^3, H^1_0(\Omega_\ep)^3}}\,,\quad \hbox{for any }\varphi\in H^1_0(Q_\varepsilon)^3\,,
\end{equation}
and compute the right hand side of (\ref{F}) by using  in (\ref{form_var})$_1$, which  gives
\begin{equation}\label{equality_duality}
\begin{array}{l}
\displaystyle
\left\langle F_{\varepsilon},\varphi\right\rangle_{H^{-1}(Q_\varepsilon)^3, H^1_0(Q_\ep)^3}=\displaystyle
-\int_{\Omega_\ep}D u_\ep: D\mathcal{R}^\varepsilon\varphi\,dx+2N^2\int_{\Omega_\ep}{\rm rot}(w_\ep)\cdot \mathcal{R}^\varepsilon\varphi\,dx+\int_{\Omega_\ep}f'\cdot (\mathcal{R}^\varepsilon\varphi)'\,dx \,.
\end{array}\end{equation}
Using Lemma \ref{Lem:estimates_velocity} for fixed $\ep$, we see that it is a bounded functional on $H^1_0(Q_\ep)$ (see Step  2 of the proof), and in fact $F_\ep\in H^{-1}(Q_\ep)^3$. Moreover, ${\rm div}(\varphi)=0$ implies $\left\langle F_{\varepsilon},\varphi\right\rangle=0\,,$ and the DeRham theorem gives the existence of $P_\varepsilon\in L^{2}_0(Q_\varepsilon)$ with $F_\varepsilon=\nabla P_\varepsilon$.
\\

{\it Step 2. Estimates for dilated and extended pressure.} Consider $\widetilde P_\ep$ obtained from $P_\ep$ by using the change of variables (\ref{dilatacion_p}). By using the Ne${\breve{\rm c}}$as inequality   for $\widetilde P_\ep\in L^{2}_0(\Omega)$, then 
$$\|\widetilde P_\ep\|_{L^{2}(\Omega)}\leq C\|\nabla_z  \widetilde P_\ep\|_{H^{-1}(\Omega)^3}\leq C\|\nabla_{h_\ep} \widetilde P_\ep\|_{H^{-1}(\Omega)^3},$$
and thus, to prove (\ref{esti_P}), it is enough to prove the second estimate  in (\ref{esti_P}) for $\nabla_{h_\ep} \widetilde P_\ep$. \\

Let us prove it. For any $\widetilde\varphi\in H^1_0(\Omega)^3$, using the change of variables (\ref{dilatacion_p}), we have
$$\begin{array}{rl}\left\langle \nabla_{h_\ep}\widetilde P_\ep, \widetilde \varphi\right\rangle_{H^{-1}(\Omega)^3, H^1 _0(\Omega)^3}=&\displaystyle -\int_{\Omega}\widetilde P_\ep\,{\rm div}_{h_\ep}(\widetilde\varphi)\,dz\\
\noame
=&\displaystyle-h_\ep^{-1}\int_{Q_\ep}P_\ep\,{\rm div}(\varphi)\,dx=h_\ep^{-1}\left\langle \nabla  P_\ep,  \varphi\right\rangle_{H^{-1}(Q_\ep)^3, H^1_0(Q_\ep)^3}.
\end{array}$$
Then, using the identification (\ref{equality_duality}) of $F_\ep$, we get
$$
\begin{array}{l}
\displaystyle
\left\langle \nabla_{h_\ep}\widetilde P_\ep,\widetilde \varphi\right\rangle_{H^{-1}(\Omega)^3, H^1_0(\Omega)^3}=\displaystyle h_\ep^{-1}\left(
-\int_{\Omega_\ep}D u_\ep: D\mathcal{R}^{\varepsilon}\varphi\,dx+2N^2\int_{\Omega_\ep}{\rm rot}(w_\ep)\cdot \mathcal{R}^{\varepsilon}\varphi\,dx+\int_{\Omega_\ep}f'\cdot (\mathcal{R}^{\varepsilon}\varphi)'\,dx\right) \,,
\end{array}
$$
and applying the change of variables (\ref{dilatacion_p}), we get
\begin{equation}\label{equality_duality2}
\begin{array}{l}
\displaystyle
\left\langle \nabla_{h_\ep}\widetilde P_\ep,\widetilde \varphi\right\rangle_{H^{-1}(\Omega)^3, H^1_0(\Omega)^3}=\displaystyle  
-\int_{\widetilde \Omega_\ep}D_{h_\ep} \widetilde u_\ep: D_{h_\ep}  \mathcal{\widetilde R}^{\varepsilon}\widetilde \varphi\,dz+2N^2\int_{\widetilde\Omega_\ep}{\rm rot}_{h_\ep}(\widetilde w_\ep)\cdot \mathcal{\widetilde R}^{\varepsilon}\widetilde\varphi\,dz+\int_{\widetilde \Omega_\ep}f'\cdot (\mathcal{R}^{\varepsilon}\widetilde\varphi)'\,dz  \,.
\end{array}
\end{equation}
Let us now estimate the right-hand side of (\ref{equality_duality2}).
From the Cauchy-Schwarz inequality,  using estimates for $\widetilde u_\ep$ in (\ref{estim_Du_tilde}), for $\widetilde w_\ep$ in (\ref{estim_Dw_tilde}), the assumption of $f$ given in (\ref{fep}) and estimates of the rescaled restricted operator (\ref{estim_restricted2}),  we   obtain
$$
\begin{array}{rl}
\displaystyle
\left|\int_{\widetilde \Omega_\ep}D_{h_\ep} \widetilde u_\ep: D_{h_\ep}  \mathcal{\widetilde R}^{\varepsilon}\widetilde \varphi\,dz\right|\leq &\displaystyle C\|D_{h_\ep} \widetilde u_\ep\|_{L^2(\widetilde \Omega_\ep)^{3\times 3}}\|D_{h_\ep} \mathcal{\widetilde R}^\ep \widetilde \varphi\|_{L^2(\widetilde \Omega_\ep)^{3\times 3}}\leq C \|\widetilde \varphi\|_{H^1_0(\Omega)^3},\\
\noame
\displaystyle \left|
2N^2\int_{\widetilde\Omega_\ep}{\rm rot}_{h_\ep}(\widetilde w_\ep)\cdot \mathcal{\widetilde R}^{\varepsilon}\widetilde\varphi\,dz
\right|\leq &\displaystyle\|D_{h_\ep} \widetilde w_\ep\|_{L^2(\widetilde\Omega_\ep)^{3\times 3}}\| \mathcal{\widetilde R}^\ep \widetilde\varphi \|_{L^2(\widetilde \Omega_\ep)^3}\leq C \|\widetilde \varphi\|_{H^1_0(\Omega)^3},
\\
\noame
\displaystyle
 \left|\int_{\widetilde \Omega_\ep}f'\cdot (\mathcal{R}^{\varepsilon}\widetilde\varphi)'\,dz\right|\leq&\displaystyle  C\| \mathcal{\widetilde R}^\ep \widetilde\varphi \|_{L^2(\widetilde \Omega_\ep)^3}\leq C\|\widetilde \varphi\|_{H^1_0(\Omega)^3}\,,
\end{array}
$$
which together with (\ref{equality_duality2}) gives $$\left|\left\langle \nabla_{h_\ep}\widetilde P_\ep,\widetilde \varphi\right\rangle_{H^{-1}(\Omega)^3, H^1_0(\Omega)^3}\right|\leq C\|\widetilde \varphi\|_{H^1_0(\Omega)^3}.$$
This implies the second estimate given in (\ref{esti_P}), which concludes the proof.
\end{proof}

\section{Adaptation of the unfolding method}\label{sec:unfolding}
This version of the unfolding method  consists in dividing the domain $\widetilde\Omega_\ep$ into squares of horizontal length $\ep$ and vertical length $\ep/h_\ep$ (see \cite{Ciora2, Cioran-book} for the original unfolding method), and is necessary to capture the microgeometry od the domain $\widetilde\Omega_\ep$. We extend the version of the unfolding method introduced in \cite[Section 4.1]{Anguiano_SG_fissure2025} described for a 2D domain to a 3D domain. To apply it, we need the following notation: for $k\in\mathbb{Z}^3$, we define $\kappa:\mathbb{R}^3\to \mathbb{Z}^3$ by
\begin{equation}\label{kappa}
\kappa(x)=k\Longleftrightarrow x\in Y_{k,1}.
\end{equation}
Remark that $\kappa$ is well defined up to a set of zero measure in $\mathbb{R}^3$, which is given by $\cup_{k\in\mathbb{Z}^3}\partial Y_{k,1}$.  Moreover, for every $\ep,\,h_\ep>0$, we have 
$$\kappa\left({x\over \ep}\right)=k\Longleftrightarrow x\in Y_{k,\ep}\quad\hbox{which is equivalent to }\quad \kappa\left({z'\over \ep},{h_\ep z_3\over \ep}\right)=k\Longleftrightarrow z\in \widetilde Y_{k,\ep}.$$

\begin{definition} Let $\widetilde \varphi$ be in $L^{q}(\widetilde\Omega_\ep)^3$,  $1\leq q<+\infty$, and $\widetilde \psi$ be in $L^{q'}(\Omega)$, $1/q+1/q'=1$. We define the functions $\widehat \varphi_\ep\in L^q(\mathbb{R}^3\times Y_f)^3$ and $\widehat \psi_\ep\in L^{q'}(\mathbb{R}^3\times Y)$ by 
\begin{eqnarray}
\displaystyle\widehat \varphi_\ep(z,y)= \widetilde \varphi\left( {\varepsilon}\kappa'\left(\frac{z'}{{\varepsilon}}, { h_\ep z_3\over \ep } \right) + {\varepsilon}y', {{\varepsilon}\over h_\ep}\kappa_3\left(\frac{z_1}{ {\varepsilon}}, { h_\ep z_2\over  \ep   } \right) +{ {\varepsilon}\over h_\ep}y_3 \right),& \hbox{a.e. }(z, y)\in  \mathbb{R}^3 \times   Y_f,&\label{def:unfolding1}\\
\noame
\displaystyle
\widehat \psi_\ep(z,y)=  \widetilde \psi\left( {\varepsilon}\kappa'\left(\frac{z'}{{\varepsilon}}, { h_\ep z_3\over \ep } \right)+ {\varepsilon}y', {{\varepsilon}\over h_\ep}\kappa_3\left(\frac{z'}{ {\varepsilon}}, { h_\ep z_3\over  \ep   } \right) +{ {\varepsilon}\over h_\ep}y_3 \right),&  \hbox{a.e. }(z, y)\in \mathbb{R}^3\times   Y,&\label{def:unfolding2}
\end{eqnarray}
assuming $\widetilde \varphi$ (resp. $\widetilde \psi$) is extended by zero outside $\widetilde\Omega_\ep$ (resp. $\Omega$), where
 the function $\kappa$ is defined by (\ref{kappa}).
\end{definition}
\begin{remark}\label{remarkCV}
The restrictions of $\widehat \varphi_{\varepsilon}$  to $\widetilde Y_{k,{\varepsilon}}\times Y_f$  (resp.  $\widehat \psi_\ep$ to $\widetilde Y_{k,{\varepsilon}}\times Y$) does not depend on $z$, while as a function of $y$ it is obtained from $\widetilde \varphi$ (resp.  from $\widetilde \psi)$ by using the change of variables 
\begin{equation}\label{CV}
y'=\frac{z'-{\varepsilon}k'}{{\varepsilon}},\quad y_3=\frac{h_\ep z_3-{\varepsilon}k_3}{{\varepsilon}},
\end{equation}
which transforms $\widetilde Y_{f_k,{\varepsilon}}$ into $Y_f$ (resp. $\widetilde Y_{k,{\varepsilon}}$ into $Y$).
\end{remark}
Next, we give some properties of the unfolded functions.
\begin{proposition} \label{properties_um}Let $\widetilde \varphi$ be in $W^{1,q}(\widetilde\Omega_\ep)^3$, $1\leq q<+\infty$, and $\widetilde \psi$ be in $L^{q'}(\Omega)$, $1/q+1/q'=1$. Then, we have
\begin{equation*}\label{relation_norms}
\begin{array}{c}
\displaystyle
\|\widehat \varphi_\ep\|_{L^q(\mathbb{R}^3\times Y_f)^3}=\|\widetilde \varphi\|_{L^q(\widetilde \Omega_\ep)^3},\\
\noame
 \|\nabla_{y'}\widehat \varphi_\ep\|_{L^q(\mathbb{R}^3\times Y_f)^{3\times 3}}=\ep\|\nabla_{z'}\widetilde \varphi\|_{L^q(\widetilde \Omega_\ep)^{3\times 3}},\quad \displaystyle \|\partial_{y_3}\widehat \varphi_\ep\|_{L^q(\mathbb{R}^3\times Y_f)^3}={\ep\over h_\ep}\|\partial_{z_3}\widetilde \varphi\|_{L^2(\widetilde \Omega_\ep)^3},\\
\noame

\|\widehat \psi_\ep\|_{L^{q'}(\mathbb{R}^3\times Y)}=\|\widetilde \psi\|_{L^{q'}(\Omega)}.

\end{array}\end{equation*}
\end{proposition}
\begin{proof} We will only make the proof for $\widehat \varphi_\ep$. The procedure for $\widehat \psi_\varepsilon$ is similar, so we omit it. Taking into account the definition (\ref{def:unfolding1}) of $\widehat{\varphi}_{\varepsilon}$, we obtain
\begin{eqnarray*}
\int_{\mathbb{R}^3\times Y_f}\left\vert D_{y'} \widehat{\varphi}_{\varepsilon}(z, y) \right\vert^q dzdy&=&\displaystyle\sum_{k\in  \mathbb{Z}^3}\int_{\widetilde Y_{k,{\varepsilon}}}\int_{Y_f}\left\vert D_{y'} \widehat{\varphi}_{\varepsilon}(z, y) \right\vert^q dy dz\\
&=&\displaystyle\sum_{k \in \mathbb{Z}^3}\int_{\widetilde Y_{k,{\varepsilon}}}\int_{Y_f}\left\vert D_{y'} \widetilde{\varphi} ({\varepsilon}k^{\prime}+{\varepsilon}y',{\varepsilon}h_\ep^{-1}k_3+ {\varepsilon}h_\ep^{-1}y_3
) \right\vert^qdydz.
\end{eqnarray*}
We observe that $\widetilde{\varphi}$ does not depend on $z$, then we can deduce
\begin{eqnarray*}
\int_{\mathbb{R}^3\times Y_f}\left\vert D_{y'} \widehat{\varphi}_{\varepsilon}(z, y) \right\vert^q dzdy
= {{\varepsilon}^3\over h_\ep}\displaystyle\sum_{k \in \mathbb{Z}^3}\int_{Y_f}\left\vert D_{y'} \widetilde{\varphi}({\varepsilon}k'+ {\varepsilon}y',{\varepsilon}h_\ep^{-1}k_3+{\varepsilon}h_\ep^{-1}y_3
) \right\vert^qdy.
\end{eqnarray*}
By the change of variables (\ref{CV}), we obtain
$$
\int_{\mathbb{R}^3\times Y_f}\left\vert D_{y'} \widehat{\varphi}_{\varepsilon}(z, y) \right\vert^q dzdy
={\varepsilon}^q
\displaystyle\sum_{k \in  \mathbb{Z}^3}\int_{\widetilde Y_{f_k,{\varepsilon}}} \left\vert D_{z'} \widetilde{\varphi}(z) \right\vert^qdz= {\varepsilon}^q\int_{\widetilde \Omega_\ep}\left\vert D_{z'} \widetilde \varphi(z) \right\vert^qdz.
$$
Thus, we get the property for $D_{y'}\widehat \varphi_\ep$.

Similarly,  we have
\begin{eqnarray*}
\int_{\mathbb{R}^3\times Y_f}\left\vert \partial_{y_3} \widehat{\varphi}_{\varepsilon}(z, y) \right\vert^qdzdy= {{\varepsilon}^3\over h_\ep}\displaystyle\sum_{k \in \mathbb{Z}^3}\int_{Y_f}\left\vert \partial_{y_3} \widetilde{\varphi} ({\varepsilon}k'+{\varepsilon}y', {\varepsilon}h_\ep^{-1}k_3+{\varepsilon}h_\ep^{-1}y_3
)\right\vert^qdy.
\end{eqnarray*}
By the change of variables (\ref{CV})  we obtain
$$
\int_{\mathbb{R}^3\times Y_f}\left\vert \partial_{y_3} \widehat{\varphi}_{\varepsilon}(z,y) \right\vert^qdzdy 
=
{{\varepsilon}^q\over h_\ep^q}\displaystyle\sum_{k \in \mathbb{Z}^3}\int_{\widetilde Y_{f_k,\ep}}\left\vert \partial_{z_3} \widetilde{\varphi} (z)\right\vert^qdz={{\varepsilon}^q\over h_\ep^q} \int_{\widetilde\Omega_\ep}\left\vert \partial_{z_3} \widetilde{\varphi}(z)\right\vert^qdz,
$$
so the the property for $\partial_{y_3}\widehat \varphi_\ep$ is proved. Finally, reasoning analogously we deduce
\begin{eqnarray*}
\int_{\mathbb{R}^3\times Y_f}\left\vert \widehat{\varphi}_{\varepsilon}(z,y)\right\vert^qdzdy = \int_{\widetilde\Omega_\ep}\left\vert \widetilde{\varphi}(z)\right\vert^qdz,
\end{eqnarray*}
and the property for $\widehat \varphi_\ep$  holds. 

\end{proof}

\begin{definition}  Consider $q=q'=2$. We define the following unfolded unknowns:
\begin{itemize}
\item[--] The unfolded velocity $\widehat u_\ep$  from  the dilated velocity $\widetilde u_\ep$ by means of (\ref{def:unfolding1}).
\item[--] The unfolded microrotation  $\widehat w_\ep$  from  the dilated microrotation $\widetilde w_\ep$ by means of (\ref{def:unfolding1}).
\item[--] The unfolded pressure $\widehat p_\ep$ from the dilated and extended pressure $\widetilde P_\ep$ by means of (\ref{def:unfolding2}).

\end{itemize}
\end{definition}

Next, we give the corresponding estimates for the unfolding unknowns. 

\begin{lemma}\label{estimates_hat}  Then, there exists a constant $C>0$ independent of $\ep$, such that $\widehat u_\ep$, $\widehat w_\ep$   and   $\widehat P_\ep$  satisfy  
\begin{equation}\label{estim_u_hat}
 \|\widehat u_\ep\|_{L^2(\mathbb{R}^3\times Y_f)^3}\leq C\ep^2,\quad 
 \|D_{y}\widehat u_\ep\|_{L^2(\mathbb{R}^3\times Y_f)^{3\times 3}}\leq C\ep^2,
 \end{equation}
 \begin{equation}\label{estim_w_hat}
 \|\widehat w_\ep\|_{L^2(\mathbb{R}^3\times Y_f)^3}\leq C\ep ,\quad 
 \|D_{y}\widehat w_\ep\|_{L^2(\mathbb{R}^3\times Y_f)^{3\times 3}}\leq C\ep,
 \end{equation}
\begin{equation}\label{estim_P_hat}
 \|\widehat P_\ep\|_{L^{2}(\mathbb{R}^3\times Y)}\leq C.
\end{equation}
 \end{lemma}
\begin{proof} 
Estimates (\ref{estim_u_hat})-(\ref{estim_P_hat}) easily follow from Proposition \ref{properties_um} with $q=q'=2$,  estimates of velocity $\widetilde u_\ep$ and microrotation $\widetilde w_\ep$ in $\widetilde\Omega_\ep$ given in Lemma \ref{Lem:estimates_velocity}, and estimate of the extended pressure $\widetilde P_\ep$ in $\Omega$ given in Lemma \ref{Estimates_extended_lemma}.

\end{proof}
 \section{Convergences of velocity, microrotation  and pressure}\label{sec:compactness}
From now on, we denote by $(\widetilde U_\ep, \widetilde W_\ep)$ the extension by zero of $(\widetilde u_\ep, \widetilde w_\ep)$ to the whole domain $\Omega$ (the velocity is zero in the obstacles). Then, estimates given in Lemma \ref{Lem:estimates_velocity} remain valid for the extension $(\widetilde U_\ep, \widetilde W_\ep)$, and $\widetilde U_\ep$  is divergence free in $\Omega$.  Here, we obtain some compactness results concerning the behavior of the sequence $(\widetilde U_\ep, \widetilde W_\ep, \widetilde P_\ep)$ and $(\widehat u_\ep,\widehat w_\ep, \widehat P_\ep)$.

\begin{lemma}\label{lemma_compactness}
Consider the extension $\widetilde U_\ep$ in $\Omega$ of the solution $\widetilde u_\ep$ of problem (\ref{form_var_tilde}). For a subsequence of $\ep$ still denoted by $\ep$,  there exist:
\begin{itemize}

\item[--]    $\widetilde u\in L^2(\Omega)^3$, with $\widetilde u_3\equiv 0$, such that
\begin{eqnarray}
&\ep^{-2}\widetilde U_\ep \rightharpoonup \widetilde u\quad\hbox{in }L^2(\Omega)^3,&\label{conv_vel_tilde}
\end{eqnarray}
\begin{equation}\label{divxproperty}
{\rm div}_{z'}\left(\int_{0}^1\widetilde u'(z)\,dz_3\right)=0\quad \hbox{in }\omega,\quad   \left(\int_{0}^1\widetilde u'(z)\,dz_3\right)\cdot n=0\quad \hbox{on }\partial\omega.
\end{equation}
\item[--]  $\widehat u\in L^2(\mathbb{R}^3; H^1_{{\rm per}}(Y_f)^3)$, with    $\widehat u=0$ in $\Omega\times  T$ and in $(\mathbb{R}^3\setminus \Omega)\times Y_f$, such that
\begin{eqnarray}
&\ep^{-2}\widehat u_\ep\rightharpoonup \widehat u\hbox{ in }L^2(\mathbb{R}^3; H^1(Y_f)^3),\quad \ep^{-2}D_y\widehat u_\ep\rightharpoonup D_y\widehat u\hbox{ in }L^2(\mathbb{R}^3\times Y_f)^{3\times 3}.\label{conv_vel_gorro}&
\end{eqnarray}
\begin{eqnarray}
& \displaystyle {\rm div}_{y}\,\widehat u(z,y)=0\quad \hbox{in }\mathbb{R}^3\times Y_f.&\label{divyproperty}\\
\noame
&\displaystyle{\rm div}_{z'}\left(\int_0^1\int_{Y_f} \widehat u'(z,y)\, dydz_3\right)=0
\quad \hbox{in }\omega,&\label{divxproperty_hat}\\
\noame
&\displaystyle \left(\int_0^1\int_{Y_f} \widehat u'(z,y)\,dydz_3 \right)\cdot n=0\quad \hbox{on }\partial\omega.& \nonumber
\end{eqnarray}
\end{itemize}
Moreover, the following relation between $\widetilde u$ and $\widehat u$ holds 
 \begin{equation}\label{relation_u_ugorro}
 \widetilde u(z)= \int_{Y_f}\widehat u(z,y)\,dy\quad  \hbox{a.e. in  }\Omega\quad\hbox{ with }\quad \int_{Y_f}\widehat u_3(z,y)\,dy=0.
 \end{equation}
\end{lemma}
\begin{proof}We divide the proof in two parts:
\begin{itemize} 
\item[--] We start with the extended velocity $\widetilde U_\ep$. The estimates of $\widetilde U_\ep$ in $\Omega$ are given in (\ref{estim_Du_tilde}), i.e. 
$$\|\widetilde U_\ep\|_{L^2( \Omega)^{3}}\leq C\ep^2,\quad 
\|D_{z'} \widetilde U_\ep\|_{L^2(\Omega)^{3\times 2}}\leq C\ep,\quad \|\partial_{z_3} \widetilde U_\ep\|_{L^2(\Omega)^{3\times 3}}\leq C\ep h_\ep,$$
 From the first estimate, we get the existence of $\widetilde u\in L^2(\Omega)^3$ such that, up to a subsequence, it holds
\begin{equation}\label{convergencia_debil_utilde}
\ep^{-2}\widetilde U_\ep \rightharpoonup \widetilde u\hbox{ in }L^2(\Omega)^3,
\end{equation}
which implies
\begin{equation}\label{convergencia_debil_utilde_dual}
\ep^{-2}{\rm div}_{z'}(\widetilde U_{\ep}')\rightharpoonup {\rm div}_{z'}(\widetilde u')\hbox{ in } L^2(0,1;H^{-1}(\omega)^3).
\end{equation} 
Since ${\rm div}_{h_\ep}(\widetilde U_\ep)=0$ in $\Omega$, multiplying by $h_\ep \ep^{-2}$, we obtain
\begin{equation}\label{div_limit1}
h_\ep \ep^{-2}{\rm div}_{z'}(\widetilde U_{\ep}')+ \ep^{-2}\partial_{z_3} \widetilde U_{\ep,3}=0\quad\hbox{in }\Omega,
\end{equation}
which, combined with (\ref{convergencia_debil_utilde_dual}), implies that $\ep^{-2}\partial_{z_3} \widetilde U_{\ep,3}$ is bounded in $L^2(0,1;H^{-1}(\omega)^3)$ and tends to zero. From the uniqueness of the limit, we have $\partial_{z_3}\widetilde u_3=0$, and so, $\widetilde u_3$ does not depend on $z_3$. Moreover, $\widetilde U_{\ep,3}=0$ on $\Gamma_0\cup\Gamma_1$, implies that $\ep^{-2}\widetilde U_{\ep,3}$ is bonded in $H^1(0,1;H^{-1}(\omega)^3)$ and then $\widetilde u_3=0$ on $\Gamma_0\cup \Gamma_1$. This implies that $\widetilde u_3\equiv 0$. This completes the proof of (\ref{conv_vel_tilde}).

Next, by considering $\widetilde \varphi\in \mathcal{D}(\omega)$  as test function in the divergence condition ${\rm div}_{h_\ep}(\widetilde U_\ep)=0$ in $\Omega$, we get
$$\int_{\Omega}\left({\rm div}_{z'}(\widetilde U_{\ep}')\, \widetilde \varphi+h_\ep^{-1}\partial_{z_3}\widetilde U_{\ep,3}\widetilde  \varphi\right)\,dz=0,$$
which, after integration by parts and multiplication by $\ep^{-2}$, gives
$$\int_{\Omega} \ep^{-{2}}\widetilde U_{\ep}'\cdot\nabla_{z'}\widetilde \varphi\,dz=0.$$
Passing to the limit by using convergence (\ref{conv_vel_tilde}), we deduce
$$\int_{\Omega}\widetilde u'\cdot \nabla_{z'}\widetilde \varphi\,dz=0,$$
and, since $\varphi$ does not depend on $z_3$, we obtain the following divergence condition (\ref{divxproperty}).

\item[--] Now we focus on the velocity $\widehat u_\ep$. From estimates of  $\widehat u_\ep$ given in (\ref{estim_u_hat}) we have the existence of $\widehat u\in L^2(\mathbb{R}^3; H^1_{\rm per}(Y_f)^3)$ satisfying, up to a subsequence,  convergences (\ref{conv_vel_gorro}). Taking into account that $\widehat u_\ep$ vanishes on $\mathbb{R}^3\times T$, we deduce that $\widehat u$ also vanishes on  $\mathbb{R}^3\times T$.  Moreover, by construction $\widehat u_\ep$ is zero outside $\widetilde\Omega_\ep$ and so, $\widehat u$ vanishes on $(\mathbb{R}^3\setminus \Omega)\times Y_f$.

Since ${\rm div}_{h_\ep}(\widetilde u_\ep)=0$ in $\widetilde \Omega_\ep$,  by applying the change of variables (\ref{CV}) we get 
$$\ep^{-1}{\rm div}_y(\widehat u_\ep)=0\quad \hbox{in }\mathbb{R}^3\times Y_f.$$
Multiplying by $\ep^{-1}$ and passing to the limit by using convergence (\ref{conv_vel_gorro}), we deduce ${\rm div}_y(\widehat u)=0$ in $\mathbb{R}^3\times Y_f$, i.e. property (\ref{divyproperty}).
\\

It remains to prove that $\widehat u$ is periodic in $y$. This follows by passing to the limit in the equality
$$\ep^{-2}\widehat u_\ep\left(z+\ep\,{\rm e}_1, -{1\over 2}, y_2, y_3\right)=\ep^{-2}\widehat u_\ep\left(z,{1\over 2}, y_2,y_3\right),$$
which is a consequence of definition (\ref{def:unfolding1}). Passing to the limit, this shows 
$$\widehat u\left(z,-{1\over 2},y_2,y_3\right)=\widehat v\left(z,{1\over 2},y_2,y_3\right),$$
and then is proved  the periodicity of $\widehat u$ with respect to $y_1$. Similarly, it can be proved the  periodicity of $\widehat u$ with respect to $y_2$. To prove the periodicity with respect to $y_3$, we consider 
$$\ep^{-2}\widehat u_\ep\left( z+{\ep\over h_\ep}{\rm e}_3,y',-{1\over 2}\right)=\ep^{-2}\widehat u_\ep\left(z,y',{1\over 2}\right),$$
and passing to the limit and  using relation (\ref{parameters}), we have 
$$\widehat u\left(z,y',-{1\over 2}\right)=\widehat v\left(z,y',{1\over 2}\right),$$
which shows the periodicity with respect to $y_3$.

Finally, relation (\ref{relation_u_ugorro}) follows from Proposition \ref{properties_um} with $q=1$, which gives
$$\int_\Omega \widetilde u(z) \,dz=\int_{\Omega\times Y_f}\widehat u(z,y)\,dzdy=\int_{\Omega}\left(\int_{Y_f}\widehat u(z,y)\,dy\right)\,dz.$$
From relation (\ref{relation_u_ugorro}) and since $\widetilde u_3\equiv 0$, it holds that $\int_{Y_f}\widehat u_3\,dy=0$.
%
%
%
%
%
Also, relation  (\ref{relation_u_ugorro}) together with (\ref{divxproperty}) gives divergence condition (\ref{divxproperty_hat}).

\end{itemize}
\end{proof}

\begin{lemma}[Convergences of microrotation] \label{lemma_compactness_micro}
Consider the extension $\widetilde W_\ep$ in $\Omega$ of the solution $\widetilde w_\ep$ of problem (\ref{form_var_tilde}). For a subsequence of $\ep$ still denoted by $\ep$,  there exist:
\begin{itemize}

\item[--]    $\widetilde w\in L^2(\Omega)^3$, with $\widetilde w_3\equiv 0$, such that
\begin{eqnarray}
& \ep^{-1}\widetilde W_\ep \rightharpoonup \widetilde w\quad\hbox{in }L^2(\Omega)^3,&\label{conv_micro_tilde}
\end{eqnarray}
 
\item[--]  $\widehat w\in L^2(\mathbb{R}^3; H^1_{{\rm per}}(Y_f)^3)$, with $\widehat w=0$ in $\Omega\times  T$ and in $(\mathbb{R}^3\setminus \Omega)\times Y_f$, such that
\begin{eqnarray}
& \ep^{-1}\widehat w_\ep\rightharpoonup \widehat w\hbox{ in }L^2(\mathbb{R}^3; H^1(Y_f)^3),\quad  \ep^{-1}D_y\widehat w_\ep\rightharpoonup D_y\widehat w\hbox{ in }L^2(\mathbb{R}^3\times Y_f)^{3\times 3}.\label{conv_micro_gorro}&
\end{eqnarray}
\end{itemize}
Moreover, the following relation between $\widetilde w$ and $\widehat w$ holds 
 \begin{equation}\label{relation_w_ugorro}
 \widetilde w(z)= \int_{Y_f}\widehat w(z,y)\,dy\quad  \hbox{a.e. in  }\Omega.
 \end{equation}
\end{lemma}
\begin{proof}
The proof of the results for the microrotation is similar to the ones of the velocity just taking into account estimates (\ref{estim_Dw_tilde}) and (\ref{estim_w_hat}), except to prove that $\widetilde w_3\equiv0$, so we only proof this. To do it, we consider as test function $\psi_\ep(z)=(0,0,\ep^{-1}\psi_3)$ in the variational formulation (\ref{form_var_tilde})$_3$, and we get
$$\begin{array}{l}\displaystyle \ep R_c\int_{\Omega}  \nabla_{z'}\widetilde W_{\ep,3}\cdot \nabla_{z'}\widetilde\psi_3\,dz+\ep h_\ep^{-2} R_c\int_{\Omega}  \partial_{z_3}\widetilde W_{\ep,3}\, \partial_{z_3}\widetilde\psi_3\,dz 
+4N^2\ep^{-1}\int_{\Omega}  \widetilde W_{\ep,3}\,\widetilde\psi_3\,dz\\
\noame\displaystyle 
=2N^2\ep^{-1}\int_{\widetilde\Omega_\ep} {\rm Rot}_{z'}(\widetilde U_\ep')\,\widetilde \psi_3\,dz.
\end{array}$$
Passing to the limit by using convergences of $\widetilde U_\ep$ given in (\ref{convergencia_debil_utilde}),  $\widetilde W_\ep$ given in (\ref{conv_micro_tilde}) and relation (\ref{parameters}), we get
$$\begin{array}{l}\displaystyle  4N^2 \int_{\Omega}  \widetilde w_{3}\,\widetilde\psi_3\,dz=0,
\end{array}$$
which implies $\widetilde w_3\equiv 0$ a.e. in $\Omega$.

\end{proof}

\begin{lemma}\label{lemma_conv_pressure}
Consider the rescaled and extended pressure $\widetilde P_\ep$. For a subsequence of $\ep$ still denoted by $\ep$, there exists  $\widetilde p\in L^{2}_0(\omega)$ independent of $z_3$, such that 
 \begin{equation}\label{conv_pressure_sub}
\widetilde P_\ep\to \widetilde p\quad\hbox{ in }L^{2}(\Omega),\quad h_\ep^{-1}\partial_{z_3}\widetilde P_\ep\rightharpoonup 0\quad\hbox{ in }H^{-1}(\Omega),
\end{equation} 
 \begin{equation}\label{conv_pressure_gorro}
\widehat P_\ep \to \widetilde p \quad\hbox{ in }L^{2}(\mathbb{R}^3\times Y).
\end{equation}

\end{lemma}
\begin{proof} Taking into account the first estimate of the pressure in (\ref{esti_P}), we deduce that there exist $\widetilde p\in L^{2}(\Omega)$ such that, up to a subsequence, 
\begin{equation}\label{conv_pressure_sub_weak}
\widetilde P_\ep\rightharpoonup \widetilde p\quad\hbox{ in }L^{2}(\Omega).
\end{equation} 
From convergence (\ref{conv_pressure_sub_weak}) we deduce that $\partial_{z_3}\widetilde P_\ep$ also converges to $\partial_{z_3}p$ in $H^{-1}(\Omega)$. Also, from the second estimate of the pressure in (\ref{esti_P}), we can deduce that  $\partial_{z_3}\widetilde P_\ep$ converges to zero in  $H^{-1}(\Omega)$.  By the uniqueness of the limit, then  we obtain $\partial_{z_3} \widetilde p=0$ and so $\widetilde p$ is independent of $z_3$. Since $\widetilde P_\ep$ has null mean value in $\Omega$, then $\widetilde p$ has null mean value in $\omega$.

Moreover, by using $\widetilde\varphi \cdot e_3$ in (\ref{form_var_tilde})$_1$, we get 
\begin{equation}\label{form_var_tilde_proof_p}
 \begin{array}{l}
\displaystyle \langle h_\ep^{-1} \partial_{z_3} \widetilde P_\ep, \widetilde \varphi_3\rangle_{H^{-1}(\Omega),H_0^1(\Omega)}
=-\int_{\Omega} h_\ep^{-2}\partial_{z_3} \widetilde U_{\ep,3}\, \partial_{z_3}\widetilde \varphi_3\,dz+2N^2\int_{ \Omega} {\rm Rot}_{z'}(\widetilde W_\ep')\, \widetilde \varphi_3\,dz,
\end{array} 
\end{equation}
where we have used the extended pressure (see (\ref{using_extension}) in Step 1 of the proof of Theorem \ref{mainthm_porous} for more details).  By using    convergences (\ref{conv_vel_tilde}) and (\ref{conv_micro_tilde}),  and taking into account that $\widetilde u_3\equiv 0$, we deduce  
$$\begin{array}{l}
\displaystyle \langle h_\ep^{-1} \partial_{z_3} \widetilde P_\ep, \widetilde \varphi_3\rangle_{H^{-1}(\Omega),H_0^1(\Omega)}
\to -\int_{\Omega}  \partial_{z_3} \widetilde u_{3}: \partial_{z_3}\widetilde \varphi_3\,dz=0,
\end{array} 
$$  and so, convergence (\ref{conv_pressure_sub})$_2$ holds.\\

  Next, following  \cite{Allaire0} adapted to the case of a thin layer, we prove that the convergence of the pressure is in fact strong. 
 Let $\widetilde \varphi_\ep, \widetilde \varphi$ be in  $H^1_0(\Omega)^3$  such that
\begin{equation}\label{strong_p_1}
\widetilde \varphi_\ep\rightharpoonup \widetilde \varphi\quad\hbox{in }H^1_0(\Omega)^3.
\end{equation}
Then, as $\widetilde p$ only depends on $z'$, we have 
$$\begin{array}{l}
\displaystyle
\left|\langle\nabla_{z}\widetilde P_\ep,\widetilde \varphi_\ep\rangle_{H^{-1}(\Omega)^3,H^1_0(\Omega)^3}-\langle\nabla_{z}\widetilde p,\widetilde \varphi\rangle_{H^{-1}(\Omega)^3,H^1_0(\Omega)^3}\right| 
\\
\noame
 \displaystyle \leq 
\left|\langle\nabla_{z}\widetilde P_\ep,\widetilde \varphi_\ep-\widetilde\varphi\rangle_{H^{-1}(\Omega)^3, H^{1}_0(\Omega)^3}\right| +\left|\langle\nabla_{z}(\widetilde P_\ep-\widetilde p),\widetilde \varphi\rangle_{H^{-1}(\Omega)^3,H^1_0(\Omega)^3}\right|.
\end{array}$$
On the one hand, using convergence (\ref{conv_pressure_sub_weak}), we have 
$$\left|\langle\nabla_{z} (\widetilde P_\ep- \widetilde p),\widetilde \varphi\rangle_{H^{-1}(\Omega)^3, H^1_0(\Omega)^3}\right|=\left|\int_\Omega\left(\widetilde P_\ep-\widetilde p\right)\,{\rm div}_z(\widetilde \varphi)\,dz\right|\to 0,\quad \hbox{as }\ep\to 0\,.$$
On the other hand, from (\ref{equality_duality2}) and proceeding similarly to the proof of Lemma \ref{Estimates_extended_lemma} but taking into account estimates (\ref{estim_restricted_tilde}) of the restricted operator $\mathcal{\widetilde R}^\ep$, we have
$$\begin{array}{rl}
\left|\langle\nabla_{z}\widetilde P_\ep,\widetilde \varphi_\ep-\widetilde \varphi\rangle_{H^{-1}(\Omega)^3,H^1_0(\Omega)^3}\right|\leq & \left|\langle\nabla_{h_\ep}\widetilde P_\ep,\widetilde \varphi_\ep-\widetilde \varphi\rangle_{H^{-1}(\Omega)^3,H^1_0(\Omega)^3}\right|\\
\noame\displaystyle
 \leq &C\left(\|\widetilde \varphi_\ep-\widetilde \varphi\|_{L^2(\Omega)^3}+\ep\|D_{h_\varepsilon} (\widetilde \varphi_\ep-\widetilde \varphi)\|_{L^2(\Omega)^{3\times 3}}\right)
 \\
\noame\displaystyle
 \leq &C\left(\|\widetilde \varphi_\ep-\widetilde \varphi\|_{L^2(\Omega)^3}+\ep h_\ep^{-1}\|D_{z} (\widetilde \varphi_\ep-\widetilde \varphi)\|_{L^2(\Omega)^{3\times 3}}\right).
\end{array}$$
The right-hand side of the previous inequality tends to zero as $\ep\to 0$, by virtue of relation (\ref{parameters}), (\ref{strong_p_1}) and the Rellich theorem. This implies that $\nabla_{z}\widetilde P_\ep\to \nabla_z \widetilde p=(\nabla_{z'} \widetilde p,0)^t$ strongly in $H^{-1}(\Omega)^3$, which together the classical Ne${\breve{\rm c}}$as inequality implies the strong convergence of the pressure $\widetilde P_\ep$ given in (\ref{conv_pressure_sub}).  Finally,   the strong convergence of $\widehat P_\varepsilon$ given in (\ref{conv_pressure_gorro}) follows from \cite[Proposition 1.9-(ii)]{Cioran-book} and the strong convergence of $\widetilde P_\varepsilon$ given in (\ref{conv_pressure_sub}).
\end{proof}

\section{Limit model}\label{sec:main_thm}
We derive the two-pressure limit Stokes micropolar model.
 \begin{theorem}\label{mainthm_porous}The triplet of limit functions $(\widehat u,\widehat w, \widetilde p)$  given in Lemmas \ref{lemma_compactness}, \ref{lemma_compactness_micro} and \ref{lemma_conv_pressure}, is the unique solution of problem 
\begin{equation}\label{system_1_2_hat}\left\{\begin{array}{rl}
-\Delta_y \widehat u + \nabla_y\widehat \pi- 2N^2 {\rm rot}_y(\widehat w)= f(z')-\nabla_{z}\widetilde p(z') & \hbox{in }\quad  Y_f, \\
\noame\displaystyle
-R_c\Delta_y \widehat w+4N^2\widehat w- 2N^2 {\rm rot}_y(\widehat u)= g(z') & \hbox{in }\quad  Y_f, \\
\noame\displaystyle
 \widehat u=0  & \hbox{ on }  T, \hbox{ for a.e. } z\in \Omega,\\
\noame
\displaystyle{\rm div}_y(\widehat u) =0& \hbox{in }\quad  Y_f, \\
\noame \displaystyle {\rm div}_{z'}\left(\int_0^1\int_{Y_f} \widehat u'(z,y)\, dydz_3\right)=0
&\hbox{ in }\omega,\\
\noame
\displaystyle \left(\int_0^1\int_{Y_f} \widehat u'(z,y)\,dydz_3 \right) n=0 & \hbox{ on }\partial\omega,
 \\
 \noame
 (\widehat u, \widehat w, \widehat \pi)\quad  Y-{\rm periodic},
\end{array}\right.
\end{equation}
with $g(z')=(g'(z'),0)^t$ and $f(z')-\nabla_{z}\widetilde p(z')=(f'(z')-\nabla_{z'}\widetilde p(z'),0)^t$. 
\end{theorem}
 
\begin{proof}[Proof of Theorem \ref{mainthm_porous}]  We divide the proof in two steps.

{\it Step 1. Variational formulation for $(\widehat u_\ep, \widehat w_\ep, \widehat P_\ep)$.} Let us first write the variational formulation satisfied by the functions $(\widehat u_\ep, \widehat w_\ep, \widehat P_\ep)$ in order to pass to the limit in step 2. 

According to Lemma \ref{lemma_compactness}, we consider $\widetilde \varphi_\ep(z)=\widehat \varphi(z,z'/\ep,h_\ep z_3/\ep)$,  as test function in (\ref{form_var_tilde})$_1$ where $\widehat \varphi(z,y)\in \mathcal{D}(\Omega;C^\infty_{\rm per}(Y)^3)$ with $\widehat \varphi(z,y)=0$ on $\Omega\times T$ and on  $(\mathbb{R}^2\setminus \Omega)\times Y$ (thus, $\widetilde \varphi_\ep(z)\in H^1_0(\widetilde\Omega_\ep)^3$). Then, we have
\begin{equation}\label{form_var_tilde_proof_v}
 \begin{array}{l}
\displaystyle\int_{\widetilde\Omega_\ep} D_{h_\ep} \widetilde u_\ep: D_{h_\ep}\widetilde \varphi_\ep\,dz-\int_{\widetilde\Omega_\ep} \widetilde p_\ep\,{\rm div}_{h_\ep}(\widetilde \varphi_\ep)\,dz
=2N^2\int_{\widetilde\Omega_\ep} {\rm rot}_{h_\ep}(\widetilde w_\ep)\cdot \widetilde \varphi_\ep\,dz +\int_{\widetilde\Omega_\ep}  f'\cdot \widetilde\varphi_\ep'\,dz\end{array} 
\end{equation}
Taking into account the extension of the pressure, we get
\begin{equation}\label{using_extension}\langle\nabla_{h_\ep} \widetilde p_\ep,  \widetilde \varphi_\ep\rangle_{H^{-1}(\widetilde\Omega_\ep)^3, H^1_0(\widetilde\Omega_\ep)^3}  =\langle\nabla_{h_\ep} \widetilde P_\ep,  \widetilde \varphi_\ep\rangle_{H^{-1}(\Omega)^3, H^1_0(\Omega)^3} =- \int_{\Omega}  \widetilde P_\ep\,{\rm div}_{h_\ep}(\widetilde \varphi_\ep)\,dz,
\end{equation}
and then,  (\ref{form_var_tilde_proof_v})  reads
\begin{equation}\label{form_var_tilde_proof_v1}
\displaystyle\int_{\widetilde\Omega_\ep} D_{h_\ep} \widetilde u_\ep: D_{h_\ep}\widetilde \varphi_\ep\,dz-\int_{\Omega} \widetilde P_\ep\,{\rm div}_{h_\ep}(\widetilde \varphi_\ep)\,dz
=2N^2\int_{\widetilde\Omega_\ep} {\rm rot}_{h_\ep}(\widetilde w_\ep)\cdot \widetilde \varphi_\ep\,dz +\int_{\widetilde\Omega_\ep}  f'\cdot \widetilde\varphi_\ep'\,dz.
\end{equation}
Taking into account the definition of $\widetilde\varphi_\ep$, it holds
\begin{equation}\label{derivadas}
\left.\begin{array}{l}
\displaystyle \partial_{z_i}\widetilde \varphi_\ep(z)=\partial_{z_i}\widehat \varphi(z,z'/\ep,h_\ep z_3/\ep)=\partial_{z_i}\widehat \varphi+\ep^{-1}\partial_{y_i}\widehat \varphi,\ i=1,2,\\
\noame
\displaystyle \partial_{z_3}\widetilde \varphi_\ep(z)=\partial_{z_3}\widehat \varphi(z,z'/\ep,h_\ep z_3/\ep)=\partial_{z_3}\widehat \varphi+h_\ep \ep^{-1}\partial_{y_3}\widehat \varphi,
\end{array}\right\}\Rightarrow \left\{\begin{array}{l}\displaystyle D_{h_\ep}\widetilde\varphi_\ep(z)=D_{h_\ep}\widehat \varphi+\ep^{-1}D_{y}\widehat \varphi,\\
\noame
\displaystyle{\rm div}_{h_\ep}(\widetilde\varphi_\ep)={\rm div}_{h_\ep}(\widehat \varphi)+\ep^{-1}{\rm div}_y(\widehat\varphi),
\end{array}\right.
\end{equation}
and then, we have that (\ref{form_var_tilde_proof_v1}) reads as follows
\begin{equation}\label{form_var_tilde_proof1}
\begin{array}{l}
\displaystyle\int_{\widetilde\Omega_\ep} D_{h_\ep} \widetilde u_\ep: D_{h_\ep}\widehat \varphi\,dz+\ep^{-1}\int_{\widetilde\Omega_\ep} D_{h_\ep} \widetilde u_\ep: D_{y}\widehat \varphi\,dz-\int_{\Omega} \widetilde P_\ep\,{\rm div}_{h_\ep}(\widehat \varphi)\,dz-\ep^{-1}\int_{\Omega} \widetilde P_\ep\,{\rm div}_{y}(\widehat \varphi)\,dz
\\
\noame
\displaystyle=2N^2\int_{\widetilde\Omega_\ep} {\rm rot}_{h_\ep}(\widetilde w_\ep)\cdot \widehat \varphi \,dz +\int_{\widetilde\Omega_\ep}  f'\cdot \widehat \varphi '\,dz.
\end{array}
\end{equation}Applying the Cauchy-Schwarz inequality and taking into account estimates (\ref{estim_Du_tilde}) and $\ep\ll h_\ep$ given in (\ref{parameters}), we get
$$\begin{array}{l}
\displaystyle 
\left|  \int_{\widetilde\Omega_\ep} D_{h_\ep} \widetilde u_\ep: D_{h_\ep}\widehat \varphi\,dz\right|\leq   C\ep\|D_{h_\ep}\widehat \varphi\|_{L^{2}(\widetilde\Omega_\ep)^{3\times 3}}\leq C\ep h_\ep^{-1}\|D_{z}\widehat \varphi\|_{L^{2}(\widetilde\Omega_\ep)^{3\times 3}}\leq C\ep h^{-1}_\ep\to 0,
\end{array}$$
and by using (\ref{conv_pressure_sub})$_2$, we have
$$\begin{array}{l}\displaystyle\int_{\Omega} \widetilde P_\ep\,{\rm div}_{h_\ep}(\widehat \varphi)\,dz=\displaystyle \int_{\Omega} \widetilde P_\ep\,{\rm div}_{z'}(\widehat \varphi')\,dz+\int_{\Omega} h_\ep^{-1}\widetilde P_\ep\,\partial_{z_3}\widehat \varphi_3\,dz =\displaystyle \int_{\Omega} \widetilde P_\ep\,{\rm div}_{z'}(\widehat \varphi)\,dz+O_\ep.
\end{array}$$
Thus, by the change of variables given in Remark \ref{remarkCV}, we obtain 
\begin{equation}\label{form_var_hat_u}
\begin{array}{l}
\displaystyle\int_{\Omega\times Y_f} \ep^{-2} D_{y} \widehat u_\ep: D_{y}\widehat \varphi\,dz-\int_{\Omega\times Y} \widehat P_\ep\,{\rm div}_{z'}(\widehat \varphi')\,dz-\ep^{-1}\int_{\Omega\times Y} \widehat P_\ep\,{\rm div}_{y}(\widehat \varphi)\,dz
\\
\noame
\displaystyle=2N^2\int_{\Omega\times Y_f} \ep^{-1} {\rm rot}_{y}(\widehat w_\ep)\cdot \widehat \varphi\,dz +\int_{\Omega\times Y_f}  f'\cdot \widehat \varphi'\,dz+O_\ep\,.
\end{array} 
\end{equation}
Next,  according to Lemma \ref{lemma_compactness_micro},  we consider $\widetilde \psi_\ep(z)= \widehat \psi(z,z'/\ep,h_\ep z_3/\ep)$,  as test function in (\ref{form_var_tilde})$_3$ where $\widehat \psi(z,y)\in \mathcal{D}(\Omega;C^\infty_{\rm per}(Y)^3)$ with $\widehat \psi(z,y)=0$ on $\Omega\times T$ and on  $(\mathbb{R}^2\setminus \Omega)\times Y$ (thus, $\widetilde \psi_\ep(z)\in H^1_0(\widetilde\Omega_\ep)^3$). Then, multiplying by $\ep^{-1}$, we have
\begin{equation}\label{form_var_tilde_w}
 \begin{array}{l}
\displaystyle \ep R_c\int_{\widetilde\Omega_\ep}  D_{h_\ep}\widetilde w_\ep: D_{h_\ep}\widetilde\psi_\ep\,dz 
+4N^2\ep^{-1}\int_{\widetilde\Omega_\ep}  \widetilde w_\ep\cdot\widetilde\psi_\ep\,dz=2N^2\ep^{-1}\int_{\widetilde\Omega_\ep} {\rm rot}_{h_\ep}(\widetilde u_\ep)\cdot\widetilde \psi_\ep\,dz + \int_{\widetilde\Omega_\ep}g'\cdot\widetilde \psi_\ep'\,dz,
\end{array} 
\end{equation}
Taking into account the definition of $\widetilde\psi_\ep$ and properties (\ref{derivadas}), we have that (\ref{form_var_tilde_w}) reads  
\begin{equation}\label{form_var_tilde_w2}
 \begin{array}{l}
\displaystyle \ep R_c\int_{\widetilde\Omega_\ep}  D_{h_\ep}\widetilde w_\ep: D_{h_\ep}\widehat \psi\,dz +  R_c\int_{\widetilde\Omega_\ep}  D_{h_\ep}\widetilde w_\ep: D_{y}\widehat\psi \,dz 
+4N^2\ep^{-1}\int_{\widetilde\Omega_\ep}  \widetilde w_\ep\cdot\widehat \psi \,dz\\
\noame
\displaystyle
=2N^2\ep^{-1}\int_{\widetilde\Omega_\ep} {\rm rot}_{h_\ep}(\widetilde u_\ep)\cdot\widehat \psi  \,dz + \int_{\widetilde\Omega_\ep}g'\cdot\widehat \psi '\,dz.
\end{array} 
\end{equation}
Applying the Cauchy-Schwarz inequality and taking into account estimates (\ref{estim_Dw_tilde}) and $\ep\ll h_\ep$ given in (\ref{parameters}), we get
$$\left|\ep R_c\int_{\widetilde\Omega_\ep}  D_{h_\ep}\widetilde w_\ep: D_{h_\ep}\widehat\psi\,dz
\right|\leq \ep C\|D_{h_\ep}\widehat \varphi\|_{L^2(\widetilde\Omega_\ep)^{3\times 3}}\leq C\ep h_\ep^{-1}\to 0,$$
and  by the change of variables given in Remark \ref{remarkCV}, we obtain 
\begin{equation}\label{form_var_tilde_w3}
 \begin{array}{l}
\displaystyle  R_c\int_{\Omega\times Y_f}  \ep^{-1} D_{y}\widehat w_\ep: D_{y}\widehat\psi\,dz 
+4N^2\ep^{-1}\int_{\Omega\times Y_f}  \widehat w_\ep\cdot\widehat\psi \,dz
\\
\noame
\displaystyle =2N^2\int_{\Omega\times Y_f} \ep^{-2}{\rm rot}_{y}(\widehat u_\ep)\cdot\widehat \psi \,dz +  \int_{\Omega\times Y_f}g'\cdot\widehat \psi '\,dz +O_\ep.
\end{array} 
\end{equation}

{\it Step 2. Passing to the limit. } Now, we   prove that the triplet of limit functions $(\widehat u, \widehat w,  \widetilde p)$   given in Lemmas \ref{lemma_compactness}, \ref{lemma_compactness_micro} and \ref{lemma_conv_pressure},  satisfies the   two pressure limit system (\ref{system_1_2_hat}), which has a unique solution $(\widehat u,\widehat w,  \widetilde p, \widehat \pi)\in L^2(\Omega;H^1_{\rm per}(Y_f)^3)\times (L^{2}_0(\omega)\cap H^1(\omega))\times L^{2}(\Omega;L^{2}_{0,{\rm per}}(Y_f))$ (see Theorem 2.4.2 and successive remark in \cite{Luka} for more detail).
\\

Let us first pass to the limit in (\ref{form_var_hat_u}). To do this, we consider $\widehat \varphi$ satisfying divergence condition ${\rm div}_y(\widehat \varphi)=0$, so (\ref{form_var_hat_u}) reads as follows
\begin{equation}\label{form_var_hat_u_pass1}
\begin{array}{l}
\displaystyle\int_{\Omega\times Y_f} \ep^{-2} D_{y} \widehat u_\ep: D_{y}\widehat \varphi\,dz-\int_{\Omega\times Y} \widehat P_\ep\,{\rm div}_{z'}(\widehat \varphi')\,dz
\\
\noame
\displaystyle=2N^2\int_{\Omega\times Y_f} \ep^{-1} {\rm rot}_{y}(\widehat w_\ep)\cdot \widehat \varphi\,dz +\int_{\Omega\times Y_f}  f'\cdot \widehat \varphi'\,dz+O_\ep\,.
\end{array} 
\end{equation}
By using convergences (\ref{conv_vel_gorro}) of $\ep^{-2}\widehat u_\ep$, (\ref{conv_micro_gorro}) of $\ep^{-1}\widehat w_\ep$ and (\ref{conv_pressure_gorro}) of $\widehat P_\ep$, passing to the limit in (\ref{form_var_hat_u_pass1}), we get
\begin{equation}\label{form_var_hat_u_pass3}
\begin{array}{l}
\displaystyle\int_{\Omega\times Y_f}  D_{y} \widehat u : D_{y}\widehat \varphi\,dz-\int_{\Omega\times Y} \widetilde p(x')\,{\rm div}_{z'}(\widehat \varphi')\,dz 
\displaystyle=2N^2\int_{\Omega\times Y_f}   {\rm rot}_{y}(\widehat w)\cdot \widehat \varphi\,dz +\int_{\Omega\times Y_f}  f'\cdot \widehat \varphi'\,dz\,,
\end{array} 
\end{equation}
which, by density, holds for every $\widehat \varphi\in \mathbb{V}$ with
$$\mathbb{V}=\left\{\begin{array}{l}
\widehat \varphi(z,y)\in L^2(\Omega;H^1_{{\rm per}}(Y)^3)\quad\hbox{ such that}\\
\noame
\displaystyle 
{\rm div}_y(\widehat \varphi)=0\quad\hbox{in }\Omega\times Y_f,\quad \widehat \varphi=0\quad\hbox{on }\Omega\times T
 \end{array}\right\}.
$$
Next, we pass to the limit in (\ref{form_var_tilde_w3}). From convergences (\ref{conv_vel_gorro}) of $\ep^{-2}\widehat u_\ep$, (\ref{conv_micro_gorro}) of $\ep^{-1}\widehat w_\ep$, we deduce that the limit equation is the following one
\begin{equation}\label{form_var_tilde_w3_limit}
 \begin{array}{l}
\displaystyle  R_c\int_{\Omega\times Y_f}  D_{y}\widehat w: D_{y}\widehat\psi\,dz 
+4N^2 \int_{\Omega\times Y_f}  \widehat w \cdot\widehat\psi \,dz  =2N^2\int_{\Omega\times Y_f}  {\rm rot}_{y}(\widehat u)\cdot\widehat \psi \,dz +  \int_{\Omega\times Y_f}g'\cdot\widehat \psi '\,dz,
\end{array} 
\end{equation}
which, by density, holds for every $\psi\in \mathbb{W} =\left\{\begin{array}{l}
\widehat \psi(z,y)\in L^2(\Omega;H^1_{{\rm per}}(Y)^3)\ \hbox{such that } \widehat \psi=0\ \hbox{on }\Omega\times T
\end{array}\right\}.
$\\

From Theorem 2.4.2 in \cite{Luka}, the variational formulation (\ref{form_var_hat_u_pass3})-(\ref{form_var_tilde_w3_limit}) admits a unique solution $(\widehat u, \widehat w)$ in $\mathbb{V}\times \mathbb{W}$. Therefore, by integration by parts, the variational formulation (\ref{form_var_hat_u_pass3})-(\ref{form_var_tilde_w3_limit}) is equivalent to the homogenized system (\ref{system_1_2_hat}).  Since $2N^2{\rm rot}_y(\widehat w)+f'\in L^2(\Omega\times Y_f)^3$ and $Y_f$ is a subset of $Y$ that is smooth and connected in the unit cell, from the regularity results for problem (\ref{system_1_2_hat}), we conclude that $\widehat u, \widehat w\in L^2(\Omega;H^2_{\rm per}(Y_f)^3)$ and then, from  (\ref{form_var_hat_u_pass3}) that  $\widetilde p\in H^1(\omega)$ (see  Part III-Theorem 2.4.2 and equations (2.4.6)-(2.4.7) in \cite{Luka},  and Lemma 4.4 in \cite{Aganovic0}).

\end{proof}

Let us define the local problems which are useful to eliminate the variable $y$ of the homogenized problem obtained in previous theorem, and then obtain a Darcy equation for the pressure $\widetilde p$.

For every $i=1,2, 3$, $k =1, 2$, we consider the following local micropolar problems:
\begin{equation}\label{local_micropolar}
\left\{
\begin{array}{ll}\displaystyle
-\Delta_y u^{i,k}+\nabla_y \pi^{i,k}-2N^2{\rm rot}_y(w^{i,k})=\widehat e_i\delta_{1k}&\hbox{in }Y_f,\\
\noame
\displaystyle
{\rm div}_y(u^{i,k})=0&\hbox{in }Y_f,\\
\noame
\displaystyle
-R_c\Delta_y w^{i,k}+4N^2w^{i,k}-2N^2{\rm rot}_y(u^{i,k})=\widehat e_i\delta_{2k}&\hbox{in }Y_f,\\
\noame
\displaystyle 
u^{i,k}=w^{i,k}=0&\hbox{in }T,\\
\noame\displaystyle 
u^{i,k}(y), w^{i,k}(y),\pi^{i,k}(y)\quad Y-{\rm periodic}.
\end{array}
\right.
\end{equation}
Here, $\widehat e_i=(\delta_{1i}, \delta_{2i}, 0)^t$, and $\delta_{ij}$  is the Kronecker delta.  It is known (see Part III-Lemma 2.5.1 in \cite{Luka}) that there exist a unique solution $(u^{i,k}, w^{i,k},\pi^{i,k})\in H^1_{{\rm per}}(Y_f)^3\times H^1_{{\rm per}}(Y_f)^3\times L^2_{0,{\rm per}}(Y_f)$ of problem (\ref{local_micropolar}), and moreover $\pi^{i,k}\in H^1(Y_f)$. One can see that $u^{3,k}=w^{3, k}=0$, and $\pi^{3,k}={\rm const}$, $k=1,2$,  since for $i = 3$ the force term is absent. 


Now, we give the main result concerning the homogenized flow.
  \begin{theorem}\label{mainthm_porous2}Consider the solution $(\widehat u,\widehat w, \widetilde p)$  of problem (\ref{system_1_2_hat}). Defining the  average velocity and microrotation by 
$$\mathcal{\widetilde U}(z')=\int_0^1\int_{Y_f}\widehat u(z,y)\,dydz_3,\quad \mathcal{\widetilde W}(z')=\int_0^1\int_{Y_f}\widehat w(z,y)\,dydz_3,$$  we have 
\begin{equation}\label{averages}
\begin{array}{l}
\displaystyle  \mathcal{\widetilde U}'(z')=K^{(1)}\left(f'(z')-\nabla_{z'}\widetilde p(z')\right)+K^{(2)}g'(z'),\quad \mathcal{\widetilde U}_{3}(z')=0\quad\hbox{in }\omega,\\
\noame
\displaystyle  \mathcal{\widetilde W}'(z')=L^{(1)}\left(f'(z')-\nabla_{z'}\widetilde p(z')\right)+L^{(2)}g'(z'),\quad \mathcal{\widetilde W}_{3}(z')=0\quad\hbox{in }\omega,
\end{array}
\end{equation}
where $K^{(k)}, L^{(k)}\in \mathbb{R}^{2\times 2}$, $k=1,2$, are matrices with coefficients
\begin{equation}\label{permeabilities}(K^{(k)})_{ij}=\int_{Y_f}u^{j,k}_i(y)\,dy,\quad (L^{(k)})_{ij}=\int_{Y_f}w^{j,k}_i(y)\,dy,\quad i,j=1,2,
\end{equation}
where $u^{j,k}, w^{j,k}$, $j,k=1, 2$, are the solutions of the local micropolar problems defined in (\ref{local_micropolar}).

Here, $\widetilde p\in H^1(\omega)\cap L^2_0(\omega)$ is the unique solution of the 2D Darcy law
\begin{equation}\label{Darcy_law}
\left\{\begin{array}{l}
\displaystyle {\rm div}_{x'}\left(K^{(1)}\left(f'(z')-\nabla_{z'}\widetilde p(z')\right)+K^{(2)}g'(z')\right)=0\quad \hbox{in }\omega,\\
\noame
\displaystyle
\left(K^{(1)}\left(f'(z')-\nabla_{z'}\widetilde p(z')\right)+K^{(2)}g'(z')\right)\cdot n=0\quad\hbox{on }\partial\omega.
\end{array}\right.
\end{equation}
\end{theorem}
\begin{proof}
We eliminate the microscopic variable $y$ in the effective problem (\ref{system_1_2_hat}). To do that, we consider the following identification:
$$\begin{array}{l}
\displaystyle \widehat u(z,y)=\sum_{j=1}^3\Big[\left(f_j(z')-\partial_{z_j}\widetilde p(z')\right)u^{j,1}(y)+g_j(z')u^{j,2}(y)\Big],
\\
\noame
\displaystyle \widehat w(z,y)=\sum_{j=1}^3\Big[\left(f_j(z')-\partial_{z_j}\widetilde p(z')\right)w^{j,1}(y)+g_j(z')w^{j,2}(y)\Big],
\\
\noame
\displaystyle \widehat \pi(z,y)=\sum_{j=1}^3\Big[\left(f_j(z')-\partial_{z_j}\widetilde p(z')\right)\pi^{j,1}(y)+g_j(z')\pi^{j,2}(y)\Big],
\end{array}$$
and thanks to the identity $ \int_{Y_f}\widehat \varphi(z,y)\,dydz_3= \widetilde \varphi(z',z_3)\,dz_3$ with $\int_{Y_f}\widehat \varphi_3\,dy=0$, satisfied by the velocity and microrotation given in Lemmas \ref{lemma_compactness} and \ref{lemma_compactness_micro}, we deduce  that 
$$\begin{array}{l}
\displaystyle  \mathcal{\widetilde U}'(z')=K^{(1)}\left(f'(z')-\nabla_{z'}\widetilde p(z')\right)+K^{(2)}g'(z'),\quad \mathcal{\widetilde U}_{3}(z')=0\quad\hbox{in }\omega,\\
\noame
\displaystyle  \mathcal{\widetilde W}'(z')=L^{(1)}\left(f'(z')-\nabla_{z'}\widetilde p(z')\right)+L^{(2)}g'(z'),\quad \mathcal{\widetilde W}_{3}(z')=0\quad\hbox{in }\omega,
\end{array}
$$
for $K^{(k)}, L^{(k)}\in \mathbb{R}^{3\times 3}, k=1, 2$. Since $u^{3,k}=w^{3,k}=0$ and $K^{(k)}, L^{(k)}$ are symmetric (see \cite{Aganovic0, Luka}), we deduce that $K^{(k)}_{3j}, L^{(k)}_{3j}=K^{(k)}_{i3}, L^{(k)}_{i3}=0$, $i,j=1, 2, 3$.   Then, we can redefine $\mathcal{\widetilde U}$ and $\mathcal{\widetilde W}$ by (\ref{averages})-(\ref{permeabilities}).

Finally, the divergence condition with respect to the variable $z'$ given in (\ref{system_1_2_hat})$_{5,6}$ and the expression of $\mathcal{\widetilde U}$ give problem (\ref{Darcy_law}), which has a unique solution since $K^{(1)}$ is positive definite (see Lemma 4.7 in \cite{Aganovic} and Part III-Theorem 2.5.2 in \cite{Luka}). Then, we can conclude that the whole sequence $(\ep^{-2}\widetilde U_\ep, \ep^{-1}\widetilde W_\ep, \widetilde P_\ep)$ converges.

\end{proof}

\subsection*{Declarations} 

\paragraph{Conflict of interest}The authors confirm that there is no conflict of interest to report.

%
%
%

\end{document}